\documentclass[]{amsart}
\usepackage{amsmath}
\usepackage{amssymb}

\theoremstyle{plain}
\newtheorem{theor}{Theorem}[section]
{}
\newtheorem{lemma}{Lemma}[theor]

\newtheorem{corol}{Corollary}[theor]
\newtheorem*{claim}{Claim}{}

\theoremstyle{definition}
\newtheorem{defin}{Definition}[section]

\newtheorem{quest}{Question}[section]

\newtheorem{condi}{Condition}[section]

\theoremstyle{remark}
\newtheorem{rem}{Remark}

\numberwithin{equation}{section}

\DeclareMathOperator{\cov}{cov}
\DeclareMathOperator{\non}{non}

\DeclareMathOperator{\dom}{domain}

\DeclareMathOperator{\hght}{HT}
\DeclareMathOperator{\sakne}{RT}
\DeclareMathOperator{\suc}{SC}
\newcommand{\IF}{\text{ if }}
\newcommand{\OR}{\text{ or }}
\newcommand{\AND}{\text{ and }}

\newcommand{\forces}[2]{\Vdash_{#1} \mbox{``} #2 \mbox{''}}

\newcommand{\Sacks}{{\mathbb S}}
\newcommand{\Reals}{{\mathbb R}}
\newcommand{\Rationals}{{\mathbb Q}}
\newcommand{\Poset}{{\mathbb P}}

\newcommand{\presup}[2]{\, ^{#1} \! #2}

\newcommand{\fomom}{\presup{\omega}{\omega}}

\newcommand{\card}[1]{\lvert #1 \rvert}

\DeclareMathOperator{\diam}{Diam}
\DeclareMathOperator{\smll}{SM}
\DeclareMathOperator{\nonmax}{Bot}
\DeclareMathOperator{\tsuc}{TS}
\DeclareMathOperator{\Top}{Top}
\newcommand{\norm}[1]{{\lVert #1 \rVert}_{\text{unif}}}
\newcommand{\SingleSeq}[1]{\vec{(#1)}}
\newcommand{\Angbr}[1]{\langle #1 \rangle}
\newcommand{\angbr}[1]{\langle\langle #1 \rangle\rangle}

\title[Unions of Rectifiable Curves]{Unions of Rectifiable Curves and
the Dimension of Banach Spaces} 
\author[J. Stepr\={a}ns]{Juris Stepr\={a}ns}
\address{Department of Mathematics, York University\\
4700 Keele Street,
North York, Ontario\\ Canada \ \ \  M3J 1P3}
\curraddr{Department of Mathematics,  University of Latvia\\
Riga, Latvia}
\email{Juris.Steprans@mathstat.yorku.ca}
\dedicatory{}
\date{}
\thanks{Research for this paper was partially supported by NSERC of Canada}
\keywords{meagre set, rectifiable curve, proper forcing}
\subjclass{Primary 03E35; Secondary 26A45}
\begin{document}
\begin{abstract}
To any metric space it is possible to associate
the cardinal invariant corresponding to the least number of
rectifiable curves in the space whose union is not meagre. It is
shown that this invariant can vary with the metric space considered,
even when restricted to the class of convex subspaces of separable
Banach spaces.  As a corollary it is obtained that it is consistent
with set theory that that any set of reals of size $\aleph_1$ is
meagre yet therer are $\aleph_1$ rectifiable curves in $\Reals^3$
whose union is not meagre. The consistency of this statement when
the phrase ``rectifiable curves'' is replaced by ``straight lines''
remains open. 
\end{abstract}
\maketitle
\tableofcontents
\bibliographystyle{plain}
\section{Introduction}
Cardinal invariants of the continuum have been studied implicitly for
almost a century and explicitly, with the help of forcing techniques,
for the past thirty years. One of the occasionally overlooked features
of the well known cardinal invariants --- such as the covering number,
additivity, cofinality and uniformity of measure and category --- is their
invariance with respect to the Polish space used to define them. This
is not true for other cardinal invariants associated with classical
real analysis. For example, the cardinal invariants associated with
covering Euclidean space with smooth surfaces depend on the
dimension of the Euclidean space being considered \cite{step.38}.
Cardinalities of maximal almost disjoint families of paths lattices of
integers also depend on the dimension of the lattice --- although in
the case results are known \cite{step.19} only in dimensions 1 and 2.
The same can be said of the problem of finding homeomorphisms of
Euclidean space which take a given $\aleph_1$-dense set to some other
$\aleph_1$-dense set \cite{step.9}. 
A striking result in this direction is due to
Cichon and Morayne \cite{cich.mora.1,mora.1,mora.2} who showed that the statement that $\aleph_n$ is a 
bound on the cardinality of the continuum is equivalent to the
assertion that there is a function $f:\Re^n \to \Re^{n+m}$ which is
onto and such that at each point at least one of the coordinate
functions is differentiable. 
The
present paper will consider another class of cardinal invariants which
occur naturally in analysis and whose values are not independent of
the metric space on which they are defined, even if the metric spaces
are taken to be separable Banach spaces.

Let ${\mathfrak Z} = (Z,d)$ be a metric space and define ${\mathcal
R}({\mathfrak Z})$ to be the $\sigma$-ideal generated by all
rectifiable curves in  ${\mathfrak Z}$. It is an easy exercise
to show that there is a $\sigma$-centred forcing which, for any $n$
covers the ground model $\Reals^n$ with countably many rectifiable
curves. Hence it is easy to find models where $\cov({\mathcal
R}({\Reals^n})) = \aleph_1$ and the continuum is as large as desired.
Theorem~\ref{t:shel}, due to S. Shelah, will establish that it is
consistent that $\non({\mathcal Null}) > \aleph_1$ yet 
 $\non({\mathcal R}(\Reals^2)) = \aleph_1$.

However, the main results of this paper will be concerned with
inequalities relating cardinal invariants of the ideal of meagre sets
with invariants of various ideals  ${\mathcal
R}({\mathfrak Z})$. The motivation for this is to be found in a
question due to P. Komjath: If $\non({\mathcal
Null})> \aleph_1$ does it follow that the union of any family of lines
in the plane of 
cardinality $\aleph_1$ is a measure zero set? The analogous  question
for meagre 
sets can also be asked: If $\non({\mathcal
Meagre})> \aleph_1$  does it follow that the union of any family of
lines in the plane of 
cardinality $\aleph_1$ is meagre?  Both of these questions remain open.
However, it will be shown that, in the second question, if one relaxes
the requirement that 
the families of
 consist of lines and require merely that they
consist of rectifiable curves (which both topologically and measure
theoretically are indistinguishable from line segments) then a
negative answer can be obtained. In other word, it will follow from
Theorem~\ref{t:main} that it is consistent that $\non({\mathcal
Meagre})> \aleph_1$ yet there is some ${\mathcal X}$ which is a family
of rectifiable curves in $\Reals^3$
such that $\card{\mathcal X} = \aleph_1$ and $\bigcup{\mathcal X}$ is
not meagre. As well, it will be shown that that the
cardinal invariant ${\mathcal
R}({\mathfrak Z})$ depends, to some extent, on the dimension of
${\mathfrak Z}$. It should be pointed out that through out this paper
by a curve is meant a one-to-one function from the unit interval to
some metric space. As will be seen in remrks in the last section, the
requirement of one-to-oneness can lead to difficulties when trying to
prove monotonicity results but it is nevertheless a natural
requirement. Not only does it eliminate Peano curves (which, of course,
are not rectifiable) but it also guarantees that the images of the curves to
be considered are all topologically equivalent to the unit interval.
This supports the assertion that the main result obtained in
Corollary~\ref{c:main} is a step
towards solving the problem of Komjath.

\section{Notation and terminology}
For the purposes of this paper, by a curve in in some metric space
${\mathfrak X} = (X,d)$ will be meant a continuous, one-to-one function
$\gamma : [0,1] \to X$. However, this notation will often be abused by
statements such as $\gamma \cap X = \emptyset$. In such cases, what
is really meant is that the image of $\gamma$ is disjoint from $X$.
Similarly, the statement that in a particular metric space there is a
family of rectifable curves whose union is not meagre really means
that the union of the images of the curves is not meagre.

For any metric space ${\mathfrak X} = (X,d)$ define a length function
on curves $\gamma$ in $X$ by  
\begin{equation}\label{n:Length}
\lambda_{\mathfrak X}(\gamma) = \sup_{0 \leq x_0 < x_1 < \ldots < x_n \leq
1}\sum_{i\in n}d(\gamma(x_i),\gamma(x_{i+1}))
\end{equation}
A curve $\gamma$ in ${\mathfrak X}$ is said to be rectifiable if
$\lambda_{\mathfrak X}(\gamma) < \infty$. 
 For any real number $\delta$ such that $0 \leq \delta < 0$
define
\begin{equation}\label{n:DeltaPaths}
{\mathbf C}({\mathfrak X},\delta) = \{A\subseteq X : \lambda_{\mathfrak X}(A) \leq \delta\}
\end{equation}
It will occasionally be useful to have a notion of length defined for sets
which are not curves by approximating the set from above by
rectifiable curves. So for an arbitrary set $A\subseteq X$ define
\begin{equation}
\lambda_{\mathfrak X}(A) =\inf\{\lambda_{\mathfrak X}(\gamma) :
\gamma\text{ is a curve and }A\subseteq \gamma\}
 \end{equation} 
The neighbourhood of the point $x \in X$ consisting of the open ball
of radius $\epsilon $ around $x$ will be denoted by
$B_{{\mathfrak X}}(x,\epsilon)$. The diametre of a subset $X$ of a
metric space will be denoted $\diam(X)$

If $\sigma$ and $\tau$ are sequences (in other words, functions whose
domains are ordinals) then $\sigma \wedge \tau$ will denote the
concatenation of $\sigma$ follows by $\tau$. 
Define $T\restriction k$ for a tree $T$.
All trees will be considered to be nonempty subtrees of some tree of the form
$\presup{\stackrel{\omega}{\smile}}{X}$ for some set $X$. If $T$ is a
tree and $t\in T$ then 
\begin{itemize}
\item the set of successors of a node $t$ in $T$ is denoted by
$\suc_T(t)$ and defined to be the set of all $x$ 
such that $t\wedge x \in T$
\item $T\angbr{t} = \{s \in T : s\subseteq t \text{ or } t \subseteq
s\}$ 
\item $T\Angbr{t} = \{s :  t \wedge s \in T\}$ 
\item the root of $T$ is denoted by $\sakne(T)$ and is defined to be
the maximal member $t$ of $T$ such that $T\angbr{t} = T$
\item the height of a tree $T$ is the supremum of the domains of its
elements and is denoted by $\hght(T)$
\item $\Top(T)$ consists of all the maximal elements of $T$ --- in
other words, $t\in \Top(T)$ if and only if $\suc_T(t) = \emptyset$
\item $\nonmax(T) = T \setminus \Top(T)$
\item $T\restriction k = \{t\in T : \card{t} \leq k\}$ for any integer $k$
\end{itemize}
When dealing with iterations of partial orders involving fusion
arguments it will be useful to have notation for trees of trees. By a
tree of trees will be understood a finite tree $\mathcal T$ such that for
each $t \in \nonmax(\mathcal T)$ the set of successors of $t$ is the
set of maximal nodes of some finite tree --- in other words,
$\suc_{\mathcal T}(t) = \Top(T)$ where $T$ is some finite tree. The
tree $T$ will be denoted by $\tsuc_{\mathcal T}(t)$.

\section{The single step partial order}
Let ${\mathfrak Y} = (Y,d)$ and ${\mathfrak Z} = (Z,e)$
be metric spaces and suppose that
${\mathbf X} = \{X_n\}_{n\in\omega}$ is a given sequence of finite
families of open subsets of $Y$. Define  $M_n({\mathbf X}) =
\prod_{i\in n}\card{X_i}$, $K_n({\mathbf X}) =
2^{nM_n({\mathbf X})}$ and 
suppose that
 $\nu_n:{\mathcal P}(X_n) \to \omega$ are functions, which are called
weak measures in \cite{shel.326}, and 
$\{\epsilon^i_j\}^{i\in\omega}_{j\in D_i}$ is an infinite triangular matrix
satisfying the following conditions:
\begin{condi}\label{c:1}
If   $A\subseteq X_n$
is such that $\nu_n(A) \geq k+1$ and
  $\bigcup \{A_i : i \in K_n({\mathbf
X})\}
= A $ then there is some $i\in K_n({\mathbf
X})$ such that
$\nu_n(A_i)\geq k$. 
\end{condi}
\begin{condi}\label{c:3}
If  $\gamma$ is a rectifiable curve $\mathfrak Y$ and
$\lambda_{\mathfrak Y}(\gamma) \leq 1$ then $\nu_n(\{a\in X_n: a\cap 
\gamma\neq \emptyset\}) =0$.   
\end{condi}
\begin{condi}\label{c:2} 
$\epsilon^i_{D_i} \leq
\epsilon^i_{D_i - 1} \leq \ldots \leq \epsilon^i_0 =
\epsilon^{i-1}_{D_{i-1}} $ and $\epsilon^i_j < 1$ for all $i$ and $j$.
\end{condi}
\begin{condi}\label{c:4} 
For any  $n\in\omega$ and any $A\subseteq X_n$ such that
$\nu_n(A) \geq k +1$ and for any function $F:A 
\to {\mathbf C}({\mathfrak Z},{\epsilon^n_{k+1}})$ there is some $B\subseteq A$
such that $\nu_n(B) \geq k $ and $\lambda_{\mathfrak Z}(\cup F[B]) \leq
\diam({\mathfrak Z})\epsilon^j_{k}$.
\end{condi}

Let $\Poset({\mathbf X},{\mathfrak Y}, {\mathfrak Z})$ be defined to
consist of all 
trees $T\subseteq \bigcup_{n\in\omega}\prod_{i\in n}X_i$ such that for
each $n\in\omega$ there is some $m\in\omega$  such that for all $t\in
T$, if $\card{t} > m$ then $\nu_{\card{t}}(\suc_T(t)) >
n(M_{\card{t}}({\mathbf X}))^{{\card{t}}^2}$.

\begin{defin}\label{d:AxA}
If $T \in \Poset({\mathbf X},{\mathfrak Y}, {\mathfrak Z})$ and
$t\in T$ define $$\nu^*_T(t) =
\nu_{\card{t}}(\suc_T(t))$$
and $$\nu^{**}_T(t) =
\frac{\nu_{T}^*(t)}{(M_{\card{t}}({\mathbf X}))^{{\card{t}}^2}}$$
and, for  $n\in \omega$,
define $\smll_n(T) = \{t\wedge\SingleSeq{x}\in T : \nu^{**}_T(t)\leq n\}$. For $T$ and $T'$ in
$\Poset({\mathbf X},{\mathfrak Y})$ define $T \leq_n T'$ if and only
if $T\subseteq T'$, $\smll_n(T')\subseteq T$ and $\nu^*_T(t) \geq n$ for
all $t \in T\setminus \smll_n(T')$.
If ${\mathcal T}$ is a tree of trees and $\sigma \in {\mathcal T}$ then define
$$\nu^*_{\mathcal T,\sigma}(t) =
\nu^*_{\tsuc_{\mathcal T}(\sigma)}(t)$$
and define
$$\nu^{**}_{\mathcal T,\sigma}(t) =
\nu^{**}_{\tsuc_{\mathcal T}(\sigma)}(t)$$
\end{defin}

\begin{lemma}\label{l:AxA}
If $T_n \in \Poset({\mathbf X},{\mathfrak Y})$ and $T_{n+1} \leq _n
T_n$ for each $n\in \omega$ then $\cap_{n\in\omega}T_n \in
\Poset({\mathbf X},{\mathfrak Y})$.
\end{lemma}
\begin{proof} This is standard. See page 365 of \cite{shel.326} or
Lemma 7.3.5 on page 340 of \cite{ba.ju.book}  for 
similar arguments. Note that $\leq_n$ corresponds to $\leq_n^*$ of
\cite{shel.326}.  
 \end{proof} 

\begin{lemma}\label{l:stand}
If $T \in \Poset({\mathbf X},{\mathfrak Y}, {\mathfrak Z})$ and $T\forces{}{\Dot{x}
\in V}$ then for all $n \in \omega$ there is some finite $S\subseteq
T$ and
$T' \leq_n T$ such
that $\Top({S})$ is a maximal antichain in $T$,
$\smll_n(T) \subseteq S \subseteq T'$ and a function $\Phi :
\Top(S) \to V$ such that $T'\angbr{s}\forces{}{\Dot{x} = \Phi(s)}$ for
each $ s\in S$.\end{lemma}
\begin{proof} Note that according to Definition~\ref{d:AxA}
$T\in\Poset({\mathbf X},{\mathfrak Y}, {\mathfrak Z})$ if and only if
$T\subseteq \prod_{i\in\omega}X_i$ is a tree and for each $n\in\omega$
there are only finitely many $t\in T$ such that $nu_T^*(t) < n$.
This allows the  standard Laver argument to be applied. See
\cite{ba.ju.book}.\end{proof} 

\begin{lemma}\label{l:cond:1}
If $S \in \Poset({\mathbf X},{\mathfrak Y}, {\mathfrak Z})$,
$\card{\sakne(T)} \geq k$ and suppose that $T$ is a finite subtree of $S$.
Suppose furthermore, that $\Psi :  \Top(T)\to K_k({\mathbf X})$. Then
there is $T' \subseteq 
T$ such that 
$\Top(T')\subseteq \Top(T)$, $\Psi$ is constant on $\Top(T')$ and 
$\nu^*_{T'}(t) \geq \nu^*_T(t) - 1$
 for all $t \in \nonmax(T')$.
\end{lemma}
\begin{proof}
Proceed by induction on the height of $T$. If $\hght(T) = 1$ then
$K_0({\mathbf X}) =1$ so that $\Psi$ is already constant to begin with.
So assume that the lemma has been established for trees of height less
than or equal to $m$ and let $T$ be a tree of height $m+1$ and let
$\Psi :\Top(T) \to M_{k}({\mathbf X})$.
Using Condition~\ref{c:1} and the fact that $K_{m}({\mathbf X})
\geq K_k({\mathbf X})$   it follows that for each $s\in
\Top(T\restriction m)$ there is $Z_s\subseteq \suc_T(s)$ such that
$\nu_m(Z_s) \geq \nu_m(Z) - 1$ and $\Psi$ is constant on $Z_s$
with value $\Psi^*(s)$. Now apply the induction hypothesis to
$T\restriction m$ and the mapping $\Psi^*$. \end{proof}

\section{Constructing the weak measures}
In this section it will be shown that for various metric spaces ${\mathfrak
Y}  = (Y,d)$ and ${\mathfrak
Z}  = (Z,d)$ there exist
sequences of weak measures ${\mathbf{X}} = \{(X_n,\nu_n)\}_{n\in\omega}$ such
that the partial order $\Poset({\mathbf{X}}, {\mathfrak
Y},{\mathfrak
Z}  = (Z,d)) $ is not empty. In order to do this, it suffices to show that for
each $D\in \omega$ and open set ${\mathcal U} \subseteq Y$ there is a
family of open sets satisfying Conditions~\ref{c:1}, \ref{c:3},
\ref{c:2} and  \ref{c:4}.

\begin{defin}\label{d:tentative}
 A be a metric space   ${\mathfrak{Y}} = (Y,d)$ will be said to be
{\em quasi-infinite dimensional}\/ if for each open set ${\mathcal U}
\subseteq Y$ there is $\delta > 0$ such that for  all $m \in\omega$
there is $A\in [{\mathcal U}]^m$ such that $d(a,b) > \delta$ for all
$\{a,b\} \in [A]^2$.
 Define a function 
$\Delta_{\mathfrak{Y}}$ from $(0,1)$ to the cardinals 
by defining $\Delta_{\mathfrak{Y}}(\epsilon) = \kappa$ if $\kappa$ is
the least cardinality of a cover of ${\mathfrak{Y}}$ by open sets of
diameter less than $\epsilon$. Recall that the metric space
$\mathfrak{Y}$ is said to be {\em totally bounded } \/ if 
$\Delta_{\mathfrak{Y}}$  takes on only integer values. Finally, define
$\mathfrak{Y}$ to be {\em reasonably geodesic} if for every $x$, $y$
and $z$ belonging to $Y$ and every $\epsilon > 0$  such that
$x$ and $y$ both belong to $B_{\mathfrak Y}(z, \epsilon)$ there is,
for each $\delta > 0$ and each set $W$ such that $\lambda_{\mathfrak
Y}(W) < \infty$, a
curve $\gamma$ such that $\{x,y\}\subseteq \gamma \subseteq
B_{\mathfrak Y}(z, \epsilon)$, $\gamma\cap W = \emptyset$ and
$\lambda_{\mathfrak Y}(\gamma) < d(x,y) + \delta$.
\end{defin}

\begin{rem} Notice that any convex subset with nonempty interior of a
Banach space of whose
dimension is at least 3 is reasonably geodesic. 
\end{rem}

\begin{lemma}\label{l:qq}
Suppose that $\delta > 0$ and $\{\gamma_i\}_{i\in k}$ are rectifiable curves
in a reasonably geodesic metric space $\mathfrak Y$ such that
$\lambda_{\mathfrak Y}(\gamma_i) < \delta$ and
such that $\gamma_i \subseteq
B_{\mathfrak Y}(x, \delta/2)$. Then there is a single rectifiable curve
$\gamma$ such that $\bigcup_{i\in k} \subseteq \gamma$ and
$\lambda_{\mathfrak Y}(\gamma) < 3k\delta $. 
\end{lemma}
\begin{proof} It is elementary to prove by induction on $k$ that if
$\delta > 0$ and $\{\gamma_i\}_{i\in k}$ are rectifiable curves
such that $\gamma_i \subseteq
B_{\mathfrak Y}(x, \delta/2)$ then, for any $\epsilon > 0$, there is a
single rectifiable curve 
$\gamma$ such that $\bigcup_{i\in k} \subseteq \gamma$ and
$\lambda_{\mathfrak Y}(\gamma) < (k-1)\delta + \sum_{i\in
k}\lambda_{\mathfrak Y}(\gamma_i)  \epsilon$. The lemma follows
immediately from this. 
\end{proof}

\begin{lemma}\label{l:comb} Let ${\mathfrak{Y}}$ be a quasi-infinite
dimensional metric space and
 ${\mathfrak{Z}}$ be a totally bounded metric space. Then for
all  open sets 
${\mathcal U} 
\subseteq Y$, $K\in\omega$, $D\in \omega$ and $\Bar{\epsilon} > 0$ there is a
sequence $\{\epsilon_i\}_{i=0}^{D_i}$ such that 
$$\epsilon_{D}\leq\epsilon_{D-1} \ldots\epsilon_{D-i} \leq
\epsilon_{D-i-1} \leq \ldots \leq \epsilon_0 \leq 
\Bar{\epsilon}$$ and   
family  ${\mathcal X}$ of open subsets of ${\mathcal Y}$ as well as a
weak measure $\nu :{\mathcal P}({\mathcal X}) \to \omega$ such that
$\nu({\mathcal X}) = D$
and $\nu$ satisfies
Conditions~\ref{c:3}, \ref{c:4} as well as the following version of
Condition~\ref{c:1}: If   $A\subseteq {\mathcal X}$
is such that $\nu(A) \geq k+1$ and
  $\bigcup \{A_i\}_{i \in K} = A $  then there is some
$i\in K$ such that
$\nu(A_i)\leq k$.
\end{lemma}
\begin{proof}
Let $\delta > 0$ be such that for every $m\in\omega$ there is
$A\in [{\mathcal U}]^m$ such that $d(a,b) > \delta$ for all $\{a,b\}
\in [A]^2$. 
Construct $\{(L_i,\epsilon_i)\}_{i=0}^{D}$ by induction on $i$. Let 
$\epsilon_0 = \Bar{\epsilon}$ and $L_0 > K/\delta + 1$. Suppose that
$\epsilon_i$ and $L_i$ have 
been constructed. 
Let $\epsilon_{i+1} =
\frac{\epsilon_i}{3L_i}$ and let $L_{i+1} > K\Delta_{\mathfrak
Z}(\epsilon_{i+1})L_i$. 

To see that this
works let $D$ be given and choose  and choose
$A\in [{\mathcal U}]^{L_D}$ such that $d(a,b) > \delta$ for all $\{a,b\}
\in [A]^2$. Then let
$\delta' > 0$ be such that  
 $d(a',b') > \delta$ for all $(a',b')$ for which there exists $\{a,b\}
\in [A]^2$ such that $d(a,a') < \delta'$ and $d(b,b') < \delta'$.
Let ${\mathcal X} = \{B_{\mathfrak Y}(a,\delta')\}_{a\in A}$.
Finally, define 
$\nu({\mathcal X}') \geq k$
if and only if  $\card{{\mathcal X}')} \geq L_{k}$ for any
${\mathcal X}'\subseteq {\mathcal X}$.

To see that the appropriate version of Condition~\ref{c:1} holds note
that if $\nu({\mathcal X}')
\geq k+1$ and ${\mathcal X}' = \cup_{i\in K}A_i$ then there is some
$i \in K$ such that $\card{A_i} \geq L_{k}$. To see that
Condition~\ref{c:3} holds let $\gamma$ be curve such that
$\lambda_{\mathfrak Y}(\gamma) \leq 1$ and 
suppose that that $\card{\{a\in A : B_{\mathfrak
Y}(a,\delta')\cap \gamma \neq \emptyset\}} \geq L_0$ or, in other words,
$\nu(\{ B_{\mathfrak
Y}(a,\delta') \in {\mathcal X} : B_{\mathfrak
Y}(a,\delta') \cap \gamma \neq \emptyset\}) \geq 1$. But then choosing $x_a
\in B_{\mathfrak
Y}(a,\delta') \cap \gamma$ for each $a \in \{a\in A : B_{\mathfrak
Y}(a,\delta')\cap \gamma \neq \emptyset\}$ yields a counterexample to the
inequality $\lambda_{\mathfrak Y}(\gamma) \leq 1$ because for any
enumeration $\{a_i\}_{i=0}^{L'}$ of $\{a\in A : B_{\mathfrak
Y}(a,\delta')\cap \gamma \neq \emptyset\}$ yields that
$$\sum_{i=0}^{L'}d(x_{a_i},x_{a_{i+1}}) \geq {L'}\delta \geq L_0\delta
\geq (K/\delta + 1)\delta > 1$$
So all that remains to be checked is Condition~\ref{c:4}.
In order to see that this holds, suppose that
$\nu({\mathcal X}') \geq
k + 1$ for some
${\mathcal X}' \subseteq {\mathcal X}$.
 Suppose that $F : {\mathcal X}' \to
{\mathbf C}({\mathfrak Z},{\epsilon_{k+1}})$ and that, without loss of
generality, 
${\mathcal X}' \in [{\mathcal X}]^{L_{k+1}}$.
 Let $\{{\mathcal O}_i\}_{i\in
\Delta_{\mathfrak Z}(\epsilon_{k+1})} $ enumerate an open cover of $Z$
with sets of diameter less than $\epsilon_{k+1}$. 
For each $B_{\mathfrak Y}(a,\delta') \in {\mathcal X}'$ choose
 $J(a)
\in \Delta_{\mathfrak Z}(\epsilon_{k+1})$ such that $F(B_{\mathcal
Y}(a,\delta'))\cap 
{\mathcal O}_{J(a)} \neq \emptyset$. Then choose $j\in
\Delta_{\mathfrak Z}(\epsilon_{k+1})$ 
such that $\card{\{B_{\mathfrak Y}(a,\delta') \in {\mathcal
X}': J(a) = j\}} \geq L_{k}$ and let
${\mathcal X}'' \in [\{B_{\mathfrak Y}(a,\delta') \in {\mathcal
X}': J(a) = j\}]^{L_k}$. Hence $\nu({\mathcal X}'') \geq k$. 
Moreover it follows from Lemma~\ref{l:qq} that
$\lambda_{\mathfrak Z}(\cup F[{\mathcal X}''])
< 3L_k\epsilon_{k+1} < \epsilon_k$.\end{proof}

\begin{corol}\label{cor:norms}
Let ${\mathfrak{Y}}$ be a separable, complete, quasi-infinite
dimensional metric space and
 ${\mathfrak{Z}} $ be a totally bounded metric space. Then there is a
matrix $\{\epsilon_i\}_{j\in D_i}^{i\in\omega}$ such 
that Conditions~\ref{c:1}, \ref{c:3}, \ref{c:2} and \ref{c:4} are
satisfied, and, furthermore, if $G \subseteq\Poset({\mathbf X},
{\mathfrak Y}, {\mathfrak Z})$ is generic over some model of set
theory $V$, then, in 
$V[G]$ there a dense $G_\delta$ subset of $\mathfrak Y$ which is
disjoint from every 
rectifiable curve in  $\mathfrak Y$ which belongs\footnote{
What is meant here, of course, is that there is a code for the dense
$G_\delta$ set and the Borel set which it codes is disjoint from every
curve with a code in the ground model. To say that a curve $\gamma$
belongs to $V$ is the same as saying that $\gamma\restriction
\Rationals \in V$.
}
 to $V$.
\end{corol}
\begin{proof}
Let $\{{\mathcal U}_n\}_{n\in\omega}$ enumerate a base of open sets for
$\mathfrak Y$.
Construct $\{\epsilon_i\}_{j\in D_i}^{i\in n}$ and $\nu_i:{\mathcal
P}(X_i)\to \omega$ by induction on $n$. Suppose that
$\{\epsilon_i\}_{j\in D_i}^{i\in n}$ and $\nu_i:{\mathcal 
P}(X_i)\to \omega$ have been constructed for $i \in n$.
Let $$D_n > n2^{(\prod_{i\in n}\card{X_i})^n}$$ and the use
Lemma~\ref{l:comb} to find a sequence a
sequence $\{\epsilon^n_i\}_{i\in D_n}$ such 
that 
$$\epsilon^{n}_{D_n}\ldots\epsilon^{n}_{i+1} \leq \epsilon^{n}_i \leq \ldots \leq \epsilon^{n}_0 \leq
\epsilon^{n-1}_{D_{n-1}}$$
as well as a family  $X_n$ of open subsets of ${\mathcal U}_n$ and a
weak measure
$\nu_n$ such that $\nu_n(X_n) = D_n$ and
Conditions~\ref{c:1}, \ref{c:3} and \ref{c:4}
are all satisfied. It follows that $\Poset({\mathbf X},
{\mathfrak Y}, {\mathfrak Z})$ is not empty.

If $G$ is $\Poset({\mathbf X},
{\mathfrak Y}, {\mathfrak Z})$ generic then let $G^*$
denote the sequence defined by $G^*(i) = s(i)$ if and only if there is
some $T\in G$ such that $s= \sakne(T)$ and $\card{s} > i$. Then let
$${\mathcal V}_G = \bigcap_{i\in \in\omega}\bigcup_{j \in
\omega\setminus i}G^*(j)$$ and note that ${\mathcal V}_G $ is a dense
$G_\delta$ in $\mathfrak Z$. To see that ${\mathcal V}_G \cap \gamma =
\emptyset$ for 
every rectifiable curve $\gamma \in V$ let $T \in \Poset({\mathbf X},
{\mathfrak Y}, {\mathfrak Z})$ and let $\gamma$ be a rectifiable curve
curve in $V$. By dividing $\gamma$ into no more than finitely many
rectifiable pieces, it may be assumed that $\lambda_{\mathfrak
Y}(\gamma) < 1$. Let $r=\card{\sakne(T)}$ and define
$$T' = \{t\in T: (\forall i > k)t(i)\cap \gamma = \emptyset\}$$
 and $T'$ is a tree. Furthermore, $T\in \Poset({\mathbf X},
{\mathfrak Y}, {\mathfrak Z})$ because, if $t\in T'$ and $\card{t} >
k$ then $\nu_n(\suc_{T'}(t)) \geq \nu_n(\suc_{T}(t)) - 1$ by
Conditions~\ref{c:1} and \ref{c:3}. It follows that
$$T\forces{\Poset({\mathbf X},
{\mathfrak Y}, {\mathfrak Z})}{\gamma \cap \bigcup_{j \in
\omega\setminus k}G^*(j) = \emptyset}$$ and hence,
$1\forces{\Poset({\mathbf X},
{\mathfrak Y}, {\mathfrak Z})}{\gamma \cap{\mathcal V}_G = \emptyset}$.
\end{proof}

It is worth noting that the hypothesis in Lemma~\ref{l:comb} as well
as  Corollary~\ref{cor:norms} can be weakened by examining the proof
of Lemma~\ref{l:comb}. Given a metric space $\mathfrak Z$, $\epsilon >
0$, integers $K$ and $D$ and $\delta > 0$ it is possible to define a
pair of sequences $\{L_i\}_{i=0}^D$ and $\{\epsilon_i\}_{i=0}^D$ as in
the proof of Lemma~\ref{l:comb} as follows: $L(K,D,\delta,\epsilon)_0
= K/\delta + 1$, 
$\epsilon(K,D,\delta,\epsilon)_0 = \epsilon$ and
$\epsilon(K,D,\delta,\epsilon)_{i+1} =
\epsilon(K,D,\delta,\epsilon)_i/3L_i$ and  
$L(K,D,\delta,\epsilon)_{i+1} = K\Delta_{\mathfrak
Z}(\epsilon(K,D,\delta,\epsilon)_{i+1})L(K,D,\delta,\epsilon)_i$. All
that was 
required of the pair of metric spaces $\mathfrak Y$ and $\mathfrak Z$
in the proof 
of Lemma~\ref{l:comb} was that in every open subset $\mathcal U$ of
$\mathfrak Y$ 
and for every $\Bar{\epsilon} > 0$  and for all 
integers $K$ and $D$ there is  $\epsilon < \Bar{\epsilon} $ and
$\delta > 0$ such that  there are $L(K,D,\delta,\epsilon)_D$ points
$\mathcal U$ which are pairwise separated by a distance of at least
$\delta$ and that $\delta L(K,D,\delta,\epsilon)_0 > 1$. To see an
example of such a pair of spaces $\mathfrak Y$ and $\mathfrak Z$ where
this occurs yet $\mathfrak Y$ is not an infinite dimensional Banach
space and $\mathfrak Z = [0,1]^3$ and construct $\mathfrak Y$ to be a
compact subspace of $\ell^p$. First choose 
blocks of independent vectors $\ell^p$ of cardinality $H_n$ all of
which have norm $1/n$. By having $H_n$ grow sufficiently quickly and
taking the convex hull of the union of all the blocks it is possible
to arrange the desired property.

\section{Trees of trees}
Trees of trees will arise in the context of of iterations of some
version of the
partial order $\Poset({\mathbf X},{\mathfrak Y}, {\mathfrak Z})$. As a
consequence, throughout this section it will be assumed that all trees
are finite subtrees of some condition in $\Poset({\mathbf
X},{\mathfrak Y}, {\mathfrak Z})$ and that whenever $\mathcal T$ is a
tree of trees and $\sigma \in \mathcal T$ then $\tsuc_{\mathcal
T}(\sigma)$ is a finite subtree of some condition in $\Poset({\mathbf
X},{\mathfrak Y}, {\mathfrak Z})$.
\begin{defin}\label{d:<}
If $\mathcal T$ and $\mathcal S$ are trees of trees then define
${\mathcal T}\prec{\mathcal S}$ if and only if
\begin{gather}
{\mathcal T}\subseteq\mathcal S\\
\Top({\mathcal T}) \subseteq\Top(\mathcal S)\\
\label{d:4.5}
(\forall \sigma\in \nonmax({\mathcal T}))
\Top(\tsuc_{{\mathcal T}}(\sigma))\subseteq\Top(\tsuc_{\mathcal
S}(\sigma) )\end{gather}\end{defin}
\begin{defin}\label{d:9}
Suppose that $\mathcal T$ and $\mathcal S$ are trees of trees and that
${\mathcal T}\prec{\mathcal S}$. For any function
$f:\omega\times\omega\times\omega \to \omega$ define 
${\mathcal T}\prec_f{\mathcal S}$ if
\begin{multline}
(\forall \sigma \in \nonmax({\mathcal T})) 
(\forall t \in \nonmax(\tsuc_{{\mathcal T}}(\sigma)))
 \nu^*_{{\mathcal T}\sigma}(t) \geq \nu^*_{{\mathcal
S}\sigma}(t) -  
f(\hght({\mathcal T}),\card{\sigma},\card{t})
\end{multline}
If $f$ is a constant function with value $k$ then $\prec_f$ will be
written as $\prec_k$. Notice that this makes sense only when ${\mathcal T} \prec \mathcal S$ because
of \eqref{d:4.5}. 
\end{defin}

When dealing with trees of trees in the arguments to follow, it will
occasionally turn out to be useful to be able to restrict the
tree which determines the successors of a given node to a smaller
level. However, it is not possible to do this directly because, for
example, it $\mathcal T$ is a tree of trees  and $t \in \tsuc_{\mathcal
T}(\emptyset)$ comes from the some level below the top then there are many
extenstions of $t$ to the top level and not all of these will agree
on how the tree of trees is to continue. In other words, $t_1$ and
$t_2$ could both be extensions of $t$ in $\tsuc_{\mathcal T}(\emptyset)$
yet ${\mathcal T}\angbr{\SingleSeq{t_1}}$ could be different than
${\mathcal T}\angbr{\SingleSeq{t_2}} $. However, if $\mathcal T$ is
carefully chosen then the number of extensions of $t$ will be much
greater than the number of possibilities for
${\mathcal T}\angbr{\SingleSeq{t'}}$ where $t'$ extends $t$. A pigeonhole
argument then makes it possible to find a large set of extensions
which all agree on the way the tree of trees should be extended ---
at least to small levels. This is the content of the next series of
technical definitions and lemmas.  

\begin{defin}\label{d:10}
 For $\sigma \in \mathcal T$ define $R_k(\sigma)$ to
be a sequence with the same domain as $\sigma$ and defined by
 $R_k(\sigma)(i) = \sigma(i)\restriction k$.
\end{defin}

\begin{defin}
For a tree of trees $\mathcal T$ define
$\hght^*({\mathcal T})$ to be the maximum of all
$\hght(\tsuc_{\mathcal T}(\sigma))$ as $\sigma$ ranges over $\mathcal
T$,
$\hght_*({\mathcal T})$ to be the minimum of all
$\hght(\tsuc_{\mathcal T}(\sigma))$ as $\sigma$ ranges over $\mathcal T$
 and define $\sakne_*({\mathcal T})$ to be the minimum of all
$\card{\sakne(\tsuc_{\mathcal T}(\sigma)}$ as $\sigma$ ranges over $\mathcal
T$. Define $\nu_*({\mathcal T})$ to be
the minumum value of
$$\{\nu_{{\mathcal T}\AND\sigma}(t) : \sigma \in {\mathcal T}, t\in
\nonmax( \tsuc_{\mathcal T}(\sigma)) \AND t\not\subseteq
\sakne(\tsuc_{\mathcal T}(\sigma))\}$$ 
\end{defin}

\begin{defin}\label{d:11}
Let $\mathcal T$ be a tree of trees and $k$ an integer. The tree of
trees $\mathcal T$ will be said to be $k$-simple if  
$\tsuc_{\mathcal T}(\sigma) \restriction k = \tsuc_{\mathcal
T}(\sigma')\restriction k$ whenever $R_k(\sigma) = R_k(\sigma')$.
\end{defin}

\begin{lemma}\label{l:30}
If  $\mathcal T$ is $k$-simple then 
$${\mathcal R}(k,\mathcal T) = \{R_k(\sigma) : \sigma \in {\mathcal T}\}$$
is a tree of trees such that $\hght^*({\mathcal R}(k,\mathcal T)) \leq k$.
\end{lemma}
\begin{proof}
For any $R_k(\sigma) \in {\mathcal R}(k,\mathcal T)$ define
$$\tsuc_{{\mathcal R}(k,\mathcal T)}(R_k( \sigma)) = \tsuc_{\mathcal
T}(\sigma)\restriction k$$ and note that, because  $\mathcal T$ is
$k$-simple, there 
is no ambiguity in this definition. The fact that  $\hght^*({\mathcal
R}(k,\mathcal T)) \leq k$ is immediate.
\end{proof}

\begin{lemma}\label{l:31}
If  $\mathcal T$ is $k$-simple, $\hght_*({\mathcal T}) \geq k$ and $\sigma 
\in {\mathcal R}(k,\mathcal T)$ then 
$${\mathcal R}(k,\mathcal T,\sigma) = \{\tau  \in {\mathcal T}:\sigma
\restriction \card{\tau} =
R_k(\tau)\}$$
is a $k$-simple tree of trees such that $\sakne_*({\mathcal R}(k,\mathcal
T,\sigma)) \geq k$.
\end{lemma}   
\begin{proof}
For any $\tau \in {\mathcal R}(k,\mathcal T,\sigma)$ define
$$\tsuc_{{\mathcal R}(k,\mathcal T,\sigma)}(\tau) =  \tsuc_{\mathcal
T}(\tau)\angbr{\sigma(\card{\tau}) }$$ It is immediate that
$\sakne^*({\mathcal R}(k,\mathcal 
T,\sigma)) \geq k$ because $\sakne(\tsuc_{\mathcal T}(\tau)) \supseteq
\sigma({\card t})$ for all $\tau \in {\mathcal R}(k,{\mathcal
T},\sigma)$ because $\hght_*({\mathcal T}) \geq k$. The $k$-simplicity
of  ${\mathcal R}(k,{\mathcal
T},\sigma)$ follows from observing that if $\tau \in {\mathcal R}(k,{\mathcal
T},\sigma)$ then $$\sakne(\tsuc_{{\mathcal R}(k,{\mathcal
T},\sigma)}(\tau)) \supseteq \sigma(\card{\tau})$$
\end{proof}

\begin{lemma}\label{l:tt}
Let $\mathcal T$ be a tree of trees and let $k$ be an integer greater
than $\hght({\mathcal T})$. It is then
possible to find a tree of trees
${\mathcal T}^k \prec_{1} \mathcal T$ 
which is $k$-simple and, furthermore, 
$$\tsuc_{{\mathcal T}^k}(\sigma)\restriction k =
\tsuc_{{\mathcal T}}(\sigma)\restriction k$$
for all $\sigma \in {\mathcal
T}^k$.
\end{lemma}
\begin{proof}
Define 
${\mathcal T}^k$ by induction on $\hght({\mathcal T})$. 
If $\hght({\mathcal T}) = 0$ then
${\mathcal T}^k$ is empty. Otherwise, let
$\hght({\mathcal T}) = n+1$ and
 define $\tsuc_{{\mathcal
T}^k}(\emptyset)\restriction k$ to be $\tsuc_{{\mathcal
T}}(\emptyset)\restriction k$. 
Now fix $t \in \Top(\tsuc_{{\mathcal T}^k}(\emptyset))$ and
use the induction hypothesis to choose $\Tilde{\mathcal T}_s \prec_1
{\mathcal T}\Angbr{\SingleSeq{s}}$ which is $k$-simple. 
Define a function $\Psi$ by letting
$\Psi(s) = {\mathcal R}(k,\Tilde{\mathcal T}_{s})$ for each $s\in
\Top(\tsuc_{\mathcal 
T}(\emptyset)\angbr{t})$. Since the number of possible
values for $\Psi(s)$ is 
not greater than 
$$(2^{\prod_{i\in \card{s}}\card{X_i}})^{\card{\mathcal T}}
\leq 2^{(n+1)M_k({\mathbf X})}\leq 2^{kM_k({\mathbf X})}$$ the
last inequality being an immediate consequence of the fact that $k
\geq \hght({\mathcal T})$. Consequently,
the number of possible values for $\Psi$ is no greater than
$K_k({\mathbf X})$ and hence, it follows from 
Lemma~\ref{l:cond:1} that there is some tree  ${\mathcal S}_t
\subseteq \tsuc_{\mathcal T}(\emptyset)\angbr{t}$ such that $\Psi$ is
constant on ${\mathcal S}_t$ and $\nu^*_{{\mathcal S}_t}(r) \geq 
 \nu^*_{{\mathcal T},
\emptyset}(r) - 1$
 for all $r \in  \nonmax({\mathcal S}_t)$.  Now define 
$$\tsuc_{{\mathcal T}^k}(\emptyset)
= \bigcup\{{\mathcal S}_t : {t \in\tsuc_{{\mathcal T}}(\emptyset)\restriction
k}\}$$ Having defined $\tsuc_{{\mathcal T}^k}(\emptyset)$, in order to
define ${\mathcal T}^k$ completely, it 
suffices to 
define ${\mathcal T}^k\Angbr{\SingleSeq{t}}) = \Psi(t)$ for each $t
\in \Top(\tsuc_{{\mathcal T}^k}(\emptyset))$ noting that the induction
hypothesis implies that ${\mathcal T}^k \prec_1 \mathcal T$. 
\end{proof}

For the next lemma recall that throughout this section, whenever a
tree of trees $\mathcal T$ is mentioned it is assumed that $\tsuc_{\mathcal 
T}(\sigma)$ is a finite subtree of some condition in $\Poset({\mathbf
X},{\mathfrak Y}, {\mathfrak Z})$ for any
$\sigma \in \mathcal T$. The integers $ M_k({\mathbf X})$ and the
reals $\epsilon^i_j$ mentioned in the statement and proof of the lemma
all refer to the parameters of $\Poset({\mathbf
X},{\mathfrak Y}, {\mathfrak Z})$.

\begin{lemma}\label{l:m}
Let $f:\omega^3 \to \omega$ be the function defined by
$f(i,j,k) = M_k({\mathbf X})^{i(i-j)}$. Suppose that 
 $\mathcal T$  is a tree of trees and
$\Phi:\Top({\mathcal T}) \to 
{Y}$ is a function. Suppose also that $\sakne_*({\mathcal T}) >
\hght({\mathcal T})$ and that
$\nu_{*}({\mathcal T}) > \hght({\mathcal T})$. It is then possible to
find a tree of trees $\Bar{\mathcal T}$ such that
$\Bar{\mathcal T} \prec_f\mathcal T$ and such that
$\lambda_{\mathfrak Y}(\Phi[\Top({\Bar{\mathcal T}})]) <
\diam({\mathfrak Y})$
\end{lemma}
\begin{proof}
Induction on the height of the tree os trees $\mathcal T$ will be used
to prove the 
following, stronger, assertion from which the lemma will follow:

\begin{claim}
Suppose that  $\mathcal T$  is a tree of trees of height $h$, 
$i = \sakne_*({\mathcal T})$, $j = \hght^*({\mathcal T})$,
$Q=\nu_*({\mathcal T})$ 
and $h < i$ amd $h < Q$.
If, additionally, $\Phi:\Top({\mathcal T}) \to 
{\mathbf C}({\mathfrak Y},{\epsilon^j_0})$ is a function then it is possible to
find a nonempty tree of trees $\Bar{\mathcal T} \prec_f \mathcal T$
such that 
$\lambda_{\mathfrak Y}(\cup\Phi[\Top({\Bar{\mathcal T}})]) <
\diam({\mathfrak Y})\epsilon^i_{Q - h}$.
\end{claim}

The lemma follows immediately from the claim upon recalling that the
$\epsilon^a_b $ 
have been chosen so that
$\epsilon^a_b < 1$
for all integers $a$ and $b$.

If the  height of ${\mathcal T}$ is zero then there is nothing to
prove so suppose that the height of $\mathcal T$ is $n + 1$ and the
claim has been established for trees of trees of height $n$ or less.
Now proceed by induction on  $\hght(\tsuc_{\mathcal
T}(\emptyset))= K$. If $K \leq i$ then $\tsuc_{\mathcal
T}(\emptyset)   $ is either empty or has only a single top node. It
is therefore possible 
to apply the induction hypothesis because either ${\mathcal T}$
 has height $n$ or $\suc_{\mathcal T}(\emptyset) = \{\SingleSeq{t}\}$ and
 so the height of ${\mathcal T}\angbr{\SingleSeq{t}}$ is essentially 					       $ n$. In particular, 
in the last case it is 
possible to apply the induction hypothesis  to ${\mathcal
T}\Angbr{\SingleSeq{t}} $ to get $\Tilde{\mathcal T} \prec_f {\mathcal
T}\Angbr{\SingleSeq{t}} $ and then note that
$f(n+1,\card{\sigma},\card{t}) \geq
f(n,\card{\sigma},\card{t})$ so that it is possible to reattach the
node $\SingleSeq{t}$ to the bottom of  $\Tilde{\mathcal T}$ to get
$\Bar{\mathcal T}$ as desired.
 
Therefore it may be assumed that  $K > i$. If $e$ is defined to be
$\hght(\tsuc_{\mathcal 
T}(\emptyset)) - 1$  then $e \geq \hght({\mathcal T})$ and so, it
is  possible to apply Lemma~\ref{l:tt}.  Hence, let ${\mathcal
T}^e$ be an $e$-simple tree such that ${\mathcal
T}^e \prec_1 {\mathcal
T}$.
Now let 
$\{\tau_w\}_{w\in L}$  enumerate
${\mathcal R}(e,{\mathcal T}^e)$ noting that  $L < M_e^{n+1}$.    
  Using  the induction
hypothesis, it is possible to construct, in $L$ steps, a sequence of
$e$-simple trees $\{{\mathcal S}_w\}_{w\in L} $ such that
\begin{gather}
\label{g:0}
{\mathcal S}_0 = {\mathcal T}^e\\
\label{g:1}
{\mathcal R}(e,{\mathcal S}_w) = {\mathcal R}(e,{\mathcal T}^e)\\
\label{g:2}
{\mathcal S}_{w+1} \prec {\mathcal S}_{w}\\
\begin{split}
\label{g:3}
{\mathcal R}(e,{\mathcal S}_{w+1}, \tau_{w+1}) \prec_{\Bar{f}}
{\mathcal R}(e,{\mathcal S}_w,\tau_{w+1})
\text{ where }\Bar{f}\\ \text{ is  
defined by } \Bar{f}(i,j,\ell) = f(i-1,j-1,\ell)
\end{split}\\
\label{g:4}
(\forall \sigma \in {\mathcal R}(e,{\mathcal S}_{w+1}))\IF
\sigma \not\subseteq \tau_{w+1} \text{ then } 
\tsuc_{{\mathcal S}_{w}}(\sigma) = 
\tsuc_{{\mathcal S}_{w+1}}(\sigma)  \\
\label{g:5}
\lambda_{\mathfrak Y}(\cup\Phi[\Top({\mathcal R}(e,{\mathcal S}_w,\tau_w))]) <
\diam({\mathfrak Y})\epsilon^{e}_{Q - (n+1)}
\end{gather}
To see why this can be done, suppose that the ${\mathcal S}_w$ have been
constructed satisfying Conditions~\eqref{g:0} to \eqref{g:5} for $
w\in u$. By the choice of $e$, either $\tsuc_{{\mathcal R}(e,{\mathcal
S}_{u-1},\tau_u)}(\emptyset)$ has height less than $e$ or $\Top(
\tsuc_{{\mathcal R}(e,{\mathcal
S}_{u-1},\tau_u)}(\emptyset) ) =
\{\tau_u(0)\wedge\SingleSeq{z}\}_{z\in Z}$ for some nonempty set $Z$.
For each $z\in Z$ the tree 
of trees ${\mathcal R}_z = {\mathcal R}(e,{\mathcal
S}_{u-1}\Angbr{\SingleSeq{z}},\tau_u)$ 
has height less than 
$n+1$ so it is possible to use the induction hypothesis to find a
subtree $\Bar{\mathcal R}_z\prec_f {\mathcal R}_z$ such that the
conditions of the claim are 
satisfied. Since $\sakne_*({\mathcal R}_z) = \Bar{e}\geq e$ it follows that
$$\lambda_{\mathfrak Y}(\cup\Phi(\Top({\Bar{\mathcal R}_z}))) <
\diam({\mathfrak Y})\epsilon^{\Bar{e}}_{Q_z - n }
$$ 
where $Q_z = \nu_*({\mathcal R}_z)$. Notice that the definition of
$\nu_*$ guarantees that $Q_z
\geq Q$ and hence, $Q_z - n > 0$.
Then define $\Phi^* : Z \to {\mathbf C}({\mathfrak
Y},{\epsilon^{\Bar{e}}_{Q_z - n }}) $ by
$\Phi^*(z) = \cup\Phi[\Top({\Bar{\mathcal R}_z})]$ and use the
Condition~\ref{c:4} for the weak measure $\nu_{\Bar{e}}$ to find $Z^*\subseteq
Z$ such that 
$\nu_{\Bar{e}}(Z^*) \geq \nu_{\Bar{e}}(Z) - 1$ and $\bigcup_{z\in
Z^*}\Phi^*(z) < 
\epsilon^{\Bar{e}}_{{Q_z - (n +1)}} \leq
\epsilon^{{e}}_{{Q_z - (n +1)}} 
$. Let ${\mathcal S}_u$ be defined 
by setting declaring $\sigma \in {\mathcal S}_u$ if and only if
$\sigma \in {\mathcal S}_{u-1}$ and, if $k \leq \card{\sigma}$  and
$ R_e(\sigma\restriction k)\subseteq \tau_u$ and
$ \sigma = \SingleSeq{z}\wedge\sigma'$ then
$\sigma'\restriction (k-1) \in \Bar{\mathcal R}_z$.
It is easy to see that ${\mathcal S}_u$ is indeed a tree of trees.
It must be checked that it is also $e$-simple and that
conditions~\eqref{g:1} to \eqref{g:5} hold. To see this suppose 
that $\sigma$ amd $\sigma'$ belong to ${\mathcal S}_u$ and that
$R_e(\sigma) = R_e(\sigma')$. 
It follows from the induction hypothesis that ${\mathcal
S}_{u-1}$ is $e$-simple and that 
$$\tsuc_{{\mathcal S}_{u-1}}(\sigma)\restriction e = 
\tsuc_{{\mathcal S}_{u-1}}(\sigma')\restriction e =
\tsuc_{{\mathcal T}^e}(\sigma)\restriction e$$
Let $t \in \Top(\tsuc_{{\mathcal S}_{u-1}}(\sigma)\restriction e)$. It
suffices to show that there is some $\tau \in {\mathcal S}_u$ such that
$\tau = \sigma \wedge \SingleSeq{t'}$ such that $t \subseteq t'$.
There are two cases to consider. First, suppose that
$R_e(\sigma)\wedge\SingleSeq{t}\not\subseteq \tau_u$. In this case it
is possible to choose  $t' \in  \tsuc_{{\mathcal T}^e}(\sigma)$
which extends $t$ and to then set $\tau = \sigma \wedge \SingleSeq{t'}$.
Otherwise, $R_e(\sigma)\wedge\SingleSeq{t}\subseteq \tau_u$. This
means that $t\subseteq \sakne(\tsuc_{\Bar{\mathcal R}_z}(\sigma)$ since
$\sakne_*(\Bar{\mathcal R}_z) > e$. It should also be observed that,
in this case, $\sigma = \SingleSeq{z}\wedge \Tilde{\sigma}$ for some
$z\in Z^*$ and $\Tilde{\sigma} \in \Bar{\mathcal R}_z $.  Hence it
suffices to choose any $\tau \in 
\suc_{\Bar{\mathcal R}_z }(\Tilde{\sigma})$ since it will automatically follow
that $R_e(\tau) = R_e(\Tilde{\sigma})\wedge\SingleSeq{t}$ and hence
$\SingleSeq{z}\wedge {\tau} \in {\mathcal S}_u$.
This same argument establishes condition~\eqref{g:4} as well.
 Condition~\eqref{g:2}
follows from the fact that ${\mathcal R}_z \prec {\mathcal
R}(e,{\mathcal S}_{u-1},\tau_u)$ for each $z \in Z^*$. In order to
verify that condition~\eqref{g:3} holds let $\sigma \in \nonmax({\mathcal
R}(e,{\mathcal S}_{u},\tau_u))$ and $t \in 
\nonmax(\tsuc_{{\mathcal
R}(e,{\mathcal S}_{u},\tau_u)}(\sigma))$. If $\sigma = \emptyset$
then either $t \subset \sakne(\tsuc_{{\mathcal
R}(e,{\mathcal S}_{u},\tau_u)}(\sigma))$, and in this case
$\nu^{*}_{{\mathcal
R}(e,{\mathcal S}_{u-1},\tau_u),\sigma}(t) = \nu^{*}_{{\mathcal
R}(e,{\mathcal S}_{u},\tau_u),\sigma}(t)$, or else $t =
\sakne(\tsuc_{{\mathcal 
R}(e,{\mathcal S}_{u},\tau_u)}(\sigma))$, and in this case
$\nu^{*}_{{\mathcal
R}(e,{\mathcal S}_{u-1},\tau_u),\sigma}(t) - 1 = \nu_e(Z)- 1\leq
\nu(Z^*) =  \nu^{*}_{{\mathcal
R}(e,{\mathcal S}_{u},\tau_u),\sigma}(t)$. On the other hand,
if $\sigma \neq \emptyset$ and $\sigma = \SingleSeq{z}\wedge \sigma'$
then $$\nu^{*}_{{\mathcal 
R}(e,{\mathcal S}_{u-1},\tau_u),\sigma}(t) = 
\nu^{*}_{{\mathcal
R}(e,{\mathcal S}_{u-1}\Angbr{\SingleSeq{z}},\tau_u),\sigma'}(t)$$
and, moreover, the choice of ${\mathcal R}_z$ guarantees that
\begin{equation}\begin{split}
\ & \ \nu^{*}_{{\mathcal
R}(e,{\mathcal S}_{u-1}\Angbr{\SingleSeq{z}},\tau_u),\sigma'}(t) - 
f(\hght({\mathcal S}_{u-1}\Angbr{\SingleSeq{z}}),\card{\sigma'},
\card{t})
\\
\leq \ &
 \nu^{*}_{{\mathcal
R}_z,\sigma'}(t) -
f(\hght({\mathcal R}(e,{\mathcal
S}_{u-1}\Angbr{\SingleSeq{z}},\tau_u)),\card{\sigma'}, 
\card{t})
\\ = \ &
\nu^{*}_{{\mathcal
R}_z,\sigma'}(t) -
f(\hght({\mathcal S}_{u}-1,\card{\sigma} -1,
\card{t}) 
\\ =  \ &
 \nu^{*}_{
{\mathcal R}(e,{\mathcal S}_{u},\tau_u)
,\sigma}(t) -
f(\hght({\mathcal R}(e,{\mathcal S}_{u},\tau_u))-1,\card{\sigma} -1,
\card{t}) 
\\ = \ &
\nu^{*}_{{\mathcal R}(e,{\mathcal S}_{u},\tau_u),\sigma}(t) -
\Bar{f}(\hght({\mathcal R}(e,{\mathcal S}_{u},\tau_u)),\card{\sigma},
\card{t})	
\end{split}\end{equation}
Condition~\eqref{g:5} follows from the fact that $Q_z \geq Q$ for each
$z\in Z^*$.

Note that ${\mathcal R}(e,{\mathcal
S}_{L}) = {\mathcal R}(e,{\mathcal T}^e) $ and that, furthermore,
for each $\sigma \in \Top({\mathcal R}(e,{\mathcal S}_L))$ and $\tau
\in {\mathcal R}(e,{\mathcal T}^e,\sigma) $ condition~\eqref{g:3} of
the induction has been applied no more than 
$M_{\card{t}}^{n+1 - \card{\tau}}$  
 times. To be precise,  if $\sigma\restriction \card{\tau}
\neq\tau_w\restriction \card{\tau}$
then  Condition~\eqref{g:4} rather than Condition~\eqref{g:3}
applies at stage $w+1$ of the construction. In other words, 
Condition~\eqref{g:3} is in force only  if $\tau \in {\mathcal
R}(e,{\mathcal T}^e,\tau_w)$  and
$\sigma\restriction \card{\tau} =\tau_w\restriction \card{\tau}$.  
By counting the number of indices $w$ such that $\sigma\restriction
\card{\tau} = \tau_w\restriction\card{\tau}$ is  possible to  conclude
that, for each $\tau \in {\mathcal 
R}(e,{\mathcal S}_L,\sigma)$ and each $t \in \tsuc_{\tau \in {\mathcal
R}(e,{\mathcal S}_L,\sigma)}(\tau)$
\begin{equation}\begin{split}\label{m:ff}
\nu^*_{{\mathcal
R}(e,{\mathcal S}_L,\sigma),\tau}(t) \geq & \nu^*_{{\mathcal
R}(e,{\mathcal T}^e,\sigma),\tau}(t) - ( M_{\card{t}}({\mathbf X})^{(n+1 -
\card{\tau})} - 1)\Bar{f}(n+1,\card{\tau},\card{t}) \\
\geq &
\nu^*_{{\mathcal
R}(e,{\mathcal T},\sigma),\tau}(t) - (M_{\card{t}}({\mathbf X})^{(n+1 -
\card{\tau})} - 1)\Bar{f}(n+1,\card{\tau},\card{t}) -
1\\
= & \nu^*_{{\mathcal
R}(e,{\mathcal T},\sigma),\tau}(t) - (M_{\card{t}}({\mathbf X})^{(n+1 -
\card{\tau})} - 1)
 M_{\card{t}}({\mathbf X})^{n(n- (\card{\tau} - 1))}
- 1\\
\geq & \nu^*_{{\mathcal
R}(e,{\mathcal T},\sigma),\tau}(t) - M_{\card{t}}({\mathbf X})^{(n+1 -
\card{\tau})}
 M_{\card{t}}({\mathbf X})^{n(n- (\card{\tau} - 1))} \\
= &
 \nu^*_{{\mathcal
R}(e,{\mathcal T},\sigma),\tau}(t) - M_{\card{t}}({\mathbf X})^{(n+1)(n+1 -
\card{\tau})}\\
= & \nu^*_{{\mathcal
R}(e,{\mathcal T},\sigma),\tau}(t) -f(\card{\mathcal
T},\card{\tau},\card{t}) 
\end{split}\end{equation} 


Because the hypothesis  of the claim and the definition of $\nu_*$
guarantee that $Q_z\geq Q > n$, 
 it follows that $\epsilon^e_{Q_z - (n+1)} \leq
\epsilon^{e}_{0}$. Hence 
there is now a natural mapping $\Psi :{\mathcal 
R}(e, {\mathcal T}^e) \to {\mathbf C}({\mathfrak Y},{\diam({\mathfrak
Y})\epsilon^{e}_{0}})$ defined by 
$\Psi(\sigma) 
=\cup\Phi[\Top({\mathcal R}(e,{\mathcal S}_{L},\sigma))]$.  Since
$\hght^*({\mathcal 
R}(e, {\mathcal T}^e)) = e$, $\sakne_*({\mathcal
R}(e, {\mathcal T}^e)) = \sakne_*({\mathcal T}) = i$
and the height of 
 $\tsuc_{{\mathcal
R}(e, {\mathcal T}^e)}(\emptyset) $ is not greater than the height of
$\tsuc_{{\mathcal 
T}}(\emptyset)$ it possible to apply the induction hypothesis again 
to ${\mathcal R}(e, {\mathcal T}^e)$ and $\Psi$.
This yields a tree ${\mathcal T}^* \prec_f {\mathcal
R}(e, {\mathcal T}^e)$ such that $\lambda_{\mathfrak
Y}(\bigcup_{\sigma \in {\mathcal T}^*}\Psi(\sigma)) <
\diam({\mathfrak Y})\epsilon^{i}_{\Bar{Q} - (n+1)}$ where $\Bar{Q} =
\nu_*({\mathcal R}(e,{\mathcal T}^e)$ but, note that $\Bar{Q} \geq Q$
so $\epsilon^{i}_{\Bar{Q} - (n+1)} \leq \epsilon^{i}_{{Q} - (n+1)}$.
 Finally, let $\Bar{\mathcal T}$ be the tree defined by
$$\Bar{\mathcal T} = \{\sigma \in {\mathcal S}_L : R_e(\sigma) \in {\mathcal T}^*\}$$
 The inequality~\eqref{m:ff} and condition~\eqref{g:1} together
guarantee that $\Bar{\mathcal T} \prec_f 
\mathcal T$ and, hence, all the requirements
of the claim are satisfied.
 \end{proof}

\section{The iteration}
Given a metric spaces ${\mathfrak Y}$ and  ${\mathfrak Z}$ and a
sequence of pairs   
${\mathbf X} = \{(X_n,\nu_n)\}_{n\in\omega}$ such for each $n$, $\nu_n$
is a weak measure
$\nu_n:{\mathcal P}(X_n) \to \omega$ satisfying Conditions~\ref{c:1}
to \ref{c:4}, let $\Poset_\alpha =
\Poset_\alpha({\mathbf X},{\mathfrak Y}, {\mathfrak Z})$ be the
countable support iteration of length $\alpha$ of the partial order 
$\Poset({\mathbf X},{\mathfrak Y},{\mathfrak Z})$. This notation will be fixed
throughout this chapter and will be used without further explanation. 
Although these  hypotheses will not be needed until
Theorem~\ref{t:main} it will also be assumed that ${\mathfrak Y}$ and
${\mathfrak Z}$ are both complete, separable and, in addition, 
${\mathfrak Y}$ is quasi-infinite dimensional and  ${\mathfrak Z}$ is
reasonably geodesic and totally bounded.

\begin{defin}\label{d:fusion}
If $\Gamma \in [\omega_2]^{< \aleph_0}$, $\psi : \Gamma \to \omega$
and $p$ and $q$ are conditions  in $\Poset_{\omega_2}$ then define $p
\leq _{\Gamma,\psi} q$ if and only if $p\restriction \gamma
\forces{\Poset_{\gamma}}{p(\gamma) \leq_{\psi(\gamma)} q(\gamma)}$ for
each $ \gamma \in \Gamma$.

If $\psi$ is a constant function with value $n$ then, as usual, $p
\leq _{\Gamma,\psi} q$ will be denoted by $p
\leq _{\Gamma,n} q$.
\end{defin}
\begin{lemma}\label{l:fusion}
If $\{p_n\}_{n\in\omega}$ is a sequence of conditions in
$\Poset_{\omega_2}$ and 
$\{\Gamma_n\}_{n\in\omega}$ is a sequence of finite subsets of
$\omega_2$ such that $\cup_{n\in\omega}\Gamma_n =
\cup_{n\in\omega}\dom{p_n}$ and $p_{n+1} \leq_{\Gamma_n,n} p_n$ then
there is some $p_{\omega} \in\Poset_{\omega_2}$ such that $p_\omega
\leq p_n$ for all $n\in\omega$.
\end{lemma}
\begin{proof} This follows from Lemma~\ref{l:AxA}
andLemma~\ref{l:stand} since Axiom~A partial 
orders satisfy the standard fusion lemma. See Definition~7.1.1 and
Lemma~7.1.3 on page 326 of \cite{ba.ju.book} \end{proof}

\begin{defin}Let $\Gamma \in [\omega_2]^{< \aleph_0}$ and let $p \in
\Poset_{{\omega_2}}$. Define $\mathcal S$ to
be good for $p$ and 
$\Gamma$ by induction on $\card{\Gamma}$. If $\card{\Gamma} = 0$ then 
$\mathcal S$ is good for $p$ and
$\Gamma$ if and only if ${\mathcal S} = \emptyset$. Now suppose that
the notion of good for $p$ and
$\Gamma$ has been defined for all $p\in \Poset_{\omega_2}$. Now, if $\Gamma
\in [\omega_2]^{n+1}$, $\zeta$ is the first element og $\Gamma$  and
$p\in \Poset_{\omega_2}$, then define 
$\mathcal S$ to be good for $p$ and 
$\Gamma$ if $\mathcal S$ is a tree of trees,  
\begin{multline}\label{good:1}
p\restriction
\zeta\forces{\Poset_{\zeta}}{\tsuc_{\mathcal
S}(\emptyset)
\subseteq p(\zeta)\AND\tsuc_{\mathcal
S}(\emptyset)
\text{is a maximal antichain in }p(\zeta)}
\end{multline}
and, for each $t \in \Top(\tsuc_{\mathcal S}(\emptyset))$,
\begin{multline}
p\restriction\zeta\forces{\Poset_{\zeta}}{p(\zeta) 
\angbr{t}\forces{\Poset}{{\mathcal S}\angbr{\SingleSeq{t}} \text{ is
good for } \Gamma\setminus \{\zeta\}  \AND p/
\Poset_{\zeta + 1}}}
\end{multline}
\end{defin}
\begin{defin}
Now, if  $\mathcal S$ is good for $p$ and 
$\Gamma$ and $\sigma \in {\mathcal S}$ then define $p[\Gamma,\sigma]$
by induction on $\hght({\mathcal S})$. If  $\hght({\mathcal S}) = 0$
then define $p[\Gamma,\sigma] = p$ so assume that  $p[\Gamma,\sigma]$ has been
defined whenever $\hght({\mathcal S}) = n$ and that $\tau \in
{\mathcal T}$ and that  $\hght({\mathcal T}) = n+1$. Let $\zeta$ be
the minimal element of $\Gamma$ and let $t \in
\tsuc({\mathcal T})$ and $\tau'$ be such that $\tau =
\SingleSeq{t}\wedge\tau'$. The definition of ${\mathcal T}$ being good for $p$
and $\Gamma$ implies that 
\begin{multline}
p\restriction\zeta
\forces{\Poset_{\zeta}}{p(\zeta)\angbr{t}\forces{\Poset}
{{\mathcal S}\angbr{\SingleSeq{t}}
\text{ is good for } \Gamma\setminus \{\zeta\}  \AND
p\restriction (\zeta, \omega_2)}}
\end{multline}
 and, hence,
\begin{equation}p\restriction\zeta\forces{\Poset_{\zeta}}{ 
p(\zeta\angbr{t})\forces{\Poset}{p/\Poset_{\zeta+1}[\Gamma\setminus
\{\zeta\},\tau'] \text{ is defined}}} 
\end{equation}
 by induction.
Therefore, it is possible to define 
\begin{equation}
p[\Gamma,\tau](\beta) =  \begin{cases}
p(\beta) & \IF \beta \in \zeta\\
p(\beta)\angbr{t}& \IF  \beta = \zeta\\
(p/\Poset_{\zeta+1})[\Gamma\setminus\{\zeta\},\tau'](\beta') & \IF \beta =
\zeta + \beta'\end{cases} 
  \end{equation}
\end{defin}

\begin{lemma}\label{l:maxac}
If $p\in \Poset_\beta$ and $\Gamma\in [\beta]^{<\aleph_0}$ and
$\mathcal S$ is a tree of trees which is good for $p$ and $\Gamma$
then for any $q \leq p$ there is some $\sigma \in \Top({\mathcal S})$
such that $q \leq p[\Gamma,\sigma]$.
\end{lemma}
\begin{proof}
Standard by induction on $\card{\Gamma}$.
\end{proof}

\begin{defin}Finally, if $p \in \Poset_{\omega_2}$ and 
$\Gamma \in [\omega_2]^{< \aleph_0}$ and $\mathcal S$ is good for $p$
and $\Gamma$ then $\mathcal S$ will be said to be $n$-sufficient for
$p$ and $\Gamma$ if 
\begin{equation}\label{e:suff}
p[\Gamma,\sigma] \restriction \alpha
\forces{\Poset_\alpha}{\tsuc_{\mathcal 
S}(\sigma) \supseteq \smll_{n + \card{{\mathcal S}\Angbr{\sigma}}}
(p(\alpha))}
\end{equation}
for each $\sigma \in \mathcal S$ and $\alpha \in \gamma$ such that
$\card{(\alpha+1)\cap \Gamma} = \card{\sigma}$.
\end{defin}
\begin{lemma}\label{l:sufffus}
Suppose that $p \in \Poset_{\omega_2}$, $\Gamma \in [\omega_2]^{< \aleph_0}$, $
k\in \omega$ and $\mathcal S$ is a tree of trees
 which is $k$-sufficient for $q$ and $\Gamma$. If $q \leq p$ is such
that $\mathcal S$ is good for $q$ and $\Gamma$ and, furthermore, 
$q[\Gamma,\sigma] \leq_{\Gamma,k} p[\Gamma,\sigma]$
for
each $\sigma \in \Top({\mathcal S})$ then $q \leq_{\Gamma,k} p$.
\end{lemma}
\begin{proof} Once again proceed by induction on $\card{\Gamma}$. In
fact, the full strength of $k$-sufficiency is not needed for  this lemma.
\end{proof}

\begin{lemma}\label{l:x}
If  $\Gamma \in [\omega_2]^{<\aleph_0}$, 
$k\in \omega$
and $p\in \Poset_{\alpha}$ is such that $p\forces{}{\Dot{x} \in V}$
then there is some $\Tilde{p} \leq_{\Gamma,k} p$ and a tree of trees
$\mathcal S$ which is $k$-sufficient for both $\Tilde{p}$ and $\Gamma$
and ${p}$ and $\Gamma$ such
such that there is some $x_\sigma$ such 
that $$\Tilde{p}[\Gamma,\sigma] \forces{}{\Dot{x} = \Check{x}_\sigma}$$ for
each $\sigma \in \Top({\mathcal S})$.
\end{lemma}
\begin{proof}
Proceed by induction on $\card{\Gamma}$ noting that if $\Gamma =
\emptyset$ then there is nothing to prove. Assuming the lemma has been
established for all ${\Gamma}$ of cardinality $ n$ let $\card{\Gamma}
= n+1$, let $\mu = \min(\Gamma)$ and $\Gamma' = \Gamma\setminus
\{\mu\}$. First  
find $q_1 \leq p \restriction \mu$ such that
$q_1 \forces{\Poset_{\mu}}{\smll_k(p(\mu)) = \Check{S}}$.
Using Lemma~\ref{l:stand} and the induction hypothesis it is possible
to find $q_2 \in \Poset_{\mu+1}$ as well as $T\supseteq S$ such that
\begin{gather}
q_2\restriction \mu 
\leq q_1\\
q_2\restriction \mu\forces{}{q_2 \leq p(\mu)}\\
q_2\forces{\Poset_{\mu}}{\Top(T) \text{ is a maximal antichain in }q_2(\mu) }\\
\begin{split}
q_2\forces{\Poset_{\mu}}{(\forall t \in \Top(\check{T}))(\exists
{\mathcal T}_t)(\exists \{x^t_\sigma\}_{\sigma \in \Top({\mathcal T}_t)})
q_2(\mu)\angbr{t}\forces{\Poset}{\check{{\mathcal T}_t}
\text{ is $k$-sufficient} \text{for } \\ \Gamma'\AND q^t}} 
 \AND
q^t \leq_{\Gamma',k} p\restriction
(\mu,\omega_2)
\AND (\forall \sigma \in \Top({\mathcal T}_t))
q^t[\Gamma,\sigma]\forces{\Poset_{\omega_2}/\Poset_{\mu +
1}}{\Check{x}_\sigma^t = \Dot{x}}\end{split}
\end{gather}

Now let $M = \sum_{t\in \Top(T)}\card{{\mathcal T}_t}$. Let $q_3
\leq q_2 \restriction \mu$ and $T'$ be such that
$$q_3\forces{\Poset_{\mu}}{T\cup \smll_M(q_2(\mu))
\subseteq \check{T'} \subseteq q_2(\mu)}$$
For each $t \in \Top({T'})$ let $\Bar{t}$ be the unique member of
$\Top(T)$ such that $\Bar{t} \subseteq t$. Define $\mathcal T$ by setting
$\tsuc_{\mathcal T}(\emptyset) = T'$ and ${\mathcal T}\Angbr{\SingleSeq{t}}
= {\mathcal T}_{\Bar{t}}$ for each $t \in \Top(T')$. If $\sigma \in
\Top({\mathcal T})$ then there are unique $t\in \Top(T')$ and $\sigma' \in
\Top({\mathcal T}_{\Bar{t}})$ such that $\sigma = \SingleSeq{t}\wedge
\sigma'$ and so define $x_\sigma = x^{\Bar{t}}_{\sigma'}$.
Let $$q(\eta) = 
\begin{cases}
q_3(\eta) & \text{if } \eta \in \mu\\
q_2(\eta) & \text{if } \eta = \mu\\
q^t(\eta) & \text{if } t\in T\AND\eta > \mu \AND
  (\exists r \in \Dot{G})\\ & r\restriction \mu
 \leq q_3\AND r\restriction \mu\forces{\Poset_\mu}{r(\mu) \leq
q_2(\mu)\angbr{t}}   
\end{cases}$$
 where $\Dot{G}$ is a name for the generic filter on $\Poset_{\omega_2}$. 
It is easy to check that $q \leq_{\Gamma,k} p$ and that $\mathcal T$
is $k$-sufficient for $q$. Moreover
$q[\Gamma,\sigma]\forces{\Poset_{\omega_2}}{\Dot{x} = x_\sigma}$ for each
$\sigma \in \Top({\mathcal T})$.
\end{proof}

\begin{corol}\label{l:suff}
Given $p \in \Poset_{\omega_2}$, $\Gamma \in [\omega_2]^{< \aleph_0}$ and $
k\in \omega$ there is $q \leq_{\Gamma, k} p$ and a tree of trees
$\mathcal S$ which is $k$-sufficient for $q$ and $\Gamma$.
\end{corol}
\begin{proof}
Simply ignore $\Dot{x}$ in Lemma~\ref{l:x}.
\end{proof}
\begin{defin}
Given $p \in \Poset_{\omega_2}$, $\Gamma \in
[\omega_2]^{<\aleph_0}$, a tree of trees ${\mathcal S}$ which is
$k$-sufficient for $p$ and $\sigma \in \Top({\mathcal S})$ define
$p \sqsubseteq_{{\mathcal S},\sigma,M} q$ if and only if
\begin{itemize}
\item $ p \leq q$
\item for each $\alpha \in \Gamma$
$$p\forces{\Poset_\alpha}{(\forall t \in
p[\Gamma,\sigma](\alpha)) \text{ either }
\nu^{**}_{p[\Gamma,\sigma]}(t) \geq M \OR
\nu^{**}_{p[\Gamma,\sigma]}(t) \geq \nu^{**}_{q[\Gamma,\sigma]}(t) -
1}$$
\item if $\tau\in \Top {\mathcal S}$ and $\tau\restriction m \neq
\sigma\restriction m$ and 
$\alpha \in \Gamma$ and $\card{\Gamma\cap (\alpha + 1)} = m$ then
$$p[\Gamma,\tau]\forces{\Poset_{\alpha}}{p[\Gamma,\tau](\alpha) =
q[\Gamma,\tau](\alpha)}$$
\end{itemize}
\end{defin}

\begin{lemma}\label{l:newkey}
Let $p \in \Poset_{\omega_2}$, $\Gamma \in [\omega_2]^{<\aleph_0}$
and suppose that ${\mathcal S}$ is $k$-sufficient for $p$ and
$\Gamma$. There is some enumeration $\{\sigma_i\}_{i\in L}$ of
$\Top({\mathcal S})$  such that
if $\{p_i\}_{i\in L}$ is a sequence of conditions
satisfying the following conditions:
\begin{itemize}
\item $p_0 = p$
\item $\mathcal S$ is good for each $p_i$ and $\Gamma$
\item $p_{i+1} \sqsubseteq_{{\mathcal S},\sigma_i,L+k} p_i$
\end{itemize}
then $p_L \leq_{\Gamma, k} p$.
\end{lemma}
\begin{proof}
Recall that $\Poset = \Poset({\mathbf X},{\mathfrak Y}, {\mathfrak Z})$
and that ${\mathbf X} = \{X_i\}_{i\in \omega}$. For each $i\in \omega$
fix a well-ordering $\leq^*$ of  $\bigcup_{i\in\omega}X_i$. Given any tree of
trees which is good for some $p\in \Poset_{\omega_2}$ the natural
enumeration of $\Top({\mathcal S})$ will be the one which corresponds
to the lexicographic ordering inherited from the ordering
$\leq^*$. It will be shown that the lemma holds if
$\{\sigma_i\}_{i\in L}$ is the natural enumeration of $\Top({\mathcal S})$. 

To establish this, proceed by induction on $\card{\Gamma}$ noting that
if $\Gamma = \emptyset$ then there is nothing to do. If the lemma is proved for
$\Gamma$ of size less $n$ let $\card{\Gamma} = n$. Let $\mu =
\min(\Gamma)$ and $\Gamma' = \Gamma \setminus \{\mu\}$ and
$\tsuc_{\mathcal S}(\emptyset) = T$. 
Notice that for each $t \in \Top(T)$ 
$$\{i \in L :  \sigma _i(0) = t\}$$
forms an interval of integers in $L$.
Therefore let $\{t_j\}_{j\in J}$ be the enumeration of $\Top(T)$
induced by $\leq^*$. 
and let $\{L(j)\}_{j\in J+1}$ be integers such that $L_0 = 0$ and
$$\{\sigma_{i}\}_{i=L_j}^{L_{j+1}-1} =\{i \in L :  \sigma _i(0) =
t_j\}$$ 

Fix $j \in J$ and let $\sigma_{j,i}\in{\mathcal
S}\Angbr{\SingleSeq{t}}$ be such that $\Angbr{\SingleSeq{t}}\wedge
\sigma_{j,i} = \sigma_{i+L_j}$. 
Therefore let $\{\sigma_{j,i}\}_{i \in L_j}$ be the
natural enumeration of 
$\{\sigma \in \Top({\mathcal S}\Angbr{\SingleSeq{t}}) :
\SingleSeq{t}\wedge\sigma \in \Top({\mathcal S})\}$.
Define  $p_{i+ L_j}[\Gamma,\SingleSeq{t}] = p_{j,i}$.
The definition of
$k$-sufficient entails that ${\mathcal S}\Angbr{\SingleSeq{t}}$ is
$k$-sufficient for each $p_{j,i}$ and $\Gamma'$.
Moreover, it is easy to see that $p_{j,i+1} \sqsubseteq_{{\mathcal
S}\Angbr{\SingleSeq{t_j}},\sigma_{j,i},L+k} p_{j,i}$ because this
simply requires 
restricting the domain of the universal quantifier in the definition 
of $\sqsubseteq_{{\mathcal
S}\Angbr{\SingleSeq{t_j}},\sigma_{j,i},L+k}$.
 From the
induction hypothesis it therefore follows that
$p_{j,L_{j+1}} \leq_{\Gamma',k} p_{j,0} p[\Gamma,\SingleSeq{t_j}]$.  

Since $\mathcal S$ is $k$-sufficient it follows that 
$$p\restriction \mu\forces{\Poset_\mu}{(\forall s \in p(\mu)\angbr{t_j})
\IF \card{s} \geq \card{t_j}\text{ then }\nu_{p(\mu)}(s) \geq L_{j+1} + k}$$
   and hence, the definition of $\sqsubseteq_{{\mathcal
S}\Angbr{\SingleSeq{t_j}},\sigma_{j,i},L+k}$ guarantees that 
$$p\restriction \mu\forces{\Poset_\mu}{(\forall s \in
p_{j,L_{j+1}}(\mu)\angbr{t_j}) 
\IF \card{s} \geq \card{t_j}\text{ then }\nu_{p_{j,L_{j+1}}(\mu)}(s) \geq  k}$$
and, moreover, if $j < k$ then $p_{L_k + i}[\Gamma,\SingleSeq{t_j}] =
p_{L_{j+1} - 1}[\Gamma,\SingleSeq{t_j}] $ fo rany $i < L_{k+1}$. It therefore
follows that $p_L$ is equivalent to the condition $p^*$ defined by
 $$p^*(\alpha) = \begin{cases}
p_L(\alpha) & \IF \alpha \in \mu\\
\bigcup_{j\in J}p_{j,L_j}(\mu) &\IF \alpha = \mu\\
p_{j,L_j}(\alpha) &\IF \alpha > \mu \AND p_{j,L_j}\restriction \mu +1
\in G
		 \end{cases}
$$ 
Then, since
${\mathcal S}$ is good for $\Gamma$ and $p_L$, it follows that
$$p_L\forces{\Poset_\mu}{ T \text{ is a
maximal antichain in }p_L(\mu)}$$ and hence
 $p^* \leq_{\Gamma,k} p$.
\end{proof}

\begin{lemma}\label{l:nk2}
Let $p \in \Poset_{\omega_2}$, $\Gamma \in [\omega_2]^{<\aleph_0}$
and suppose that ${\mathcal S}$ is good for $p$ and
$\Gamma$. Suppose further that $\sigma \in \Top({\mathcal S})$ and that
$q \leq p[\Gamma,\sigma]$ satisfies, for each $\alpha\in \Gamma$, that
\begin{equation}\label{e:c4-1}
{q}\forces{\Poset_\alpha}{(\forall t \in
{q}(\alpha)) \text{ either }
\nu^{**}_{\Tilde{q}}(t) \geq M\OR
\nu^{**}_{\Tilde{q}}(t) \geq \nu^{**}_{p[\Gamma,\sigma]}(t) -
1}\end{equation}
Then there exists some  $q' \in \Poset_{\omega_2}$  such that
\begin{itemize}
\item $q'\sqsubseteq_{{\mathcal
S},\sigma,M} p$ \item $\mathcal S$ is good for $q'$ 
and $\Gamma$ \item $q'[\Gamma,\sigma] \leq q$.\end{itemize}
\end{lemma}
\begin{proof}
Proceed by induction on $\card{\Gamma}$ noting that if $\Gamma = \emptyset$
there is nothing to do. Let $\mu = \min(\Gamma)$ and $T =
\tsuc_{\mathcal S}(\emptyset)$. Let $\sigma(0) = t$ and $\sigma =
\SingleSeq{t} \wedge \sigma'$.

 Using the induction    hypothesis it
is possible to find a $\Poset_{\mu+1}$ name $q^*$ for a condition
in $\Poset_{\omega_2}/\Poset_{\mu+1}$ such
that
\begin{itemize}
\item $1\forces{\Poset_{\mu + 1}}{q^*\sqsubseteq_{{\mathcal
S}\Angbr{\SingleSeq{t}},\sigma',M} p[\Gamma,\SingleSeq{t}]/\Poset{\mu + 1}}$ 
\item $1\forces{\Poset_{\mu + 1}}{\mathcal S\Angbr{\SingleSeq{t}} \text{ is good for }q^*
\AND\Gamma' \ }$  
\item $1\forces{\Poset_{\mu + 1}}{q^*[\Gamma,\sigma'] \leq q/\Poset_{\mu +
1}}$.
\end{itemize}
Now define
$$q'(\eta) = \begin{cases}
q(\eta) & \IF \eta\in \mu \\
(p[\Gamma,\sigma](\eta) \setminus \left(p(\eta)\Angbr{\sigma(j)})\right)
\cup {q}(\eta)
& \IF \eta = \mu \\
p(\eta) & \IF \mu\in\eta \AND q\restriction (\mu + 1) \notin G_{\mu + 1}\\
q^*(\eta) & \IF \mu\in\eta \AND q\restriction (\mu + 1) \in G_{\mu + 1}
					  \end{cases}$$
where $G_{\mu + 1}$ is a name for the generic filter on $\Poset_{\mu +
1}$. The fact that $q'$ is the desired condition follows from the
induction hypothesis.
\end{proof}

\begin{lemma}\label{new:tech}
Suppose that ${\mathcal S}$ is good for $\Gamma \in
[\omega_2]^{<\aleph_0}$ and $p\in \Poset_{\omega_2}$ and that
$\Bar{\mathcal S} \prec \mathcal S$. If
 ${p_{\Bar{\mathcal S}}}$ is defined by
$${p_{\Bar{\mathcal S}}}(\alpha) = \begin{cases}
p(\alpha) & \IF \alpha \notin \Gamma\\
\bigcup\{p(\alpha)\angbr{s} : s\in \tsuc_{\Bar{\mathcal
S}}(\tau) \AND 
& \IF \alpha \in \Gamma\\
\card{\tau} = \card{(\alpha + 1) \cap \Gamma} \AND 
{p_{\Bar{\mathcal S}}}[\Gamma,\tau]\restriction \alpha \in \Dot{G}\cap
\Poset_\alpha\} & 
		      \end{cases}$$
where $\Dot{G}$ is a name for the generic filter on
$\Poset_{\omega_2}$ then $\Bar{\mathcal S}$ is  good for $\Gamma$ and
${p_{\Bar{\mathcal S}}}$. \end{lemma}
\begin{proof} Once again, proceed by induction on $\card{\Gamma}$. 
Let $\card{\Gamma} = n$ and assume the lemma proved if $\card{\Gamma} <
n$. Let $\mu = 
\min(\Gamma)$ and $\Gamma' = \Gamma \setminus \{\mu\}$ and
$\tsuc_{\Bar{\mathcal S}}(\emptyset) = T$. For each $t \in \Top(T)$  
the induction hypothesis gives that
$\Bar{\mathcal S}\Angbr{\SingleSeq{t}}$ is good for
${p[\Gamma,\SingleSeq{t}]_{\Bar{\mathcal S}}}$ and $\Gamma'$. Since $\Top(T)
\subseteq \Top(\tsuc_{\mathcal S}(\emptyset))$ it follows that if
$q \leq {p_{\Bar{\mathcal S}}}$ then $q \leq
{p[\Gamma,\SingleSeq{t}]_{\Bar{\mathcal S}}}$ for some $t \in \Top(T)$. Hence 
$\Bar{\mathcal S}$ is  good for $\Gamma$ and
${p_{\Bar{\mathcal S}}}$.
\end{proof}

\begin{theor}\label{t:main}
If $p \in 
\Poset_\alpha = \Poset_\alpha({\mathbf
X},{\mathfrak Y}, {\mathfrak Z})$ is such that
$$p \forces{\Poset_\alpha}{{\mathcal U}\text{ is a dense }
G_\delta\text{ in }{\mathfrak Z}} $$ then there is some $q \leq 
p$ as well as a rectifiable curve $\gamma$  in ${\mathfrak Z}$ such that
$q\forces{\Poset_\alpha}{\gamma\cap
{\mathcal U}\text{ is dense in } \gamma}$.
\end{theor}
\begin{proof}  
In this theorem the hypothesis that $\mathfrak Z = (Z,e)$ is reasonably
geodesic will play an important role. For later reference, notice that
one of the consequences of this hypothesis is the following 
\begin{multline}\label{h:3} (\forall \text{ open } {\mathcal V}\subseteq
Z)(\forall m\in \omega)(\forall \delta > 0)(\exists
\{x_i\}_{i\in m})(\exists \theta > 0)
(\forall i \in m)
B_{\mathfrak Z}(x_i,\theta)\subseteq
{\mathcal V} \AND\\
(\forall (\gamma_1,\gamma_2,\ldots,\gamma_m) \in \prod_{i\in
1}^m{\mathbf C}(B_{\mathfrak Z}(x_i,\theta)))(\exists \text{ a curve }
\gamma)\lambda_{\mathfrak Z}(\gamma) < \delta \AND\gamma \supseteq
\cup_{i\in m}\gamma_i 
\end{multline}
This follows from Lemma~\ref{l:qq} by choosing $\theta < \delta/3m$
such that $B_{\mathfrak Z}(x,\theta/2)\subseteq \mathcal V$ for some
point $x$. Then, noticing that any reasonable geodesic space with more
than two points can have no isolated points, simply choose
$\{{\mathcal V}\}_{i\in m}$ to be 
disjoint nonempty open subsets of $\mathcal V$.

A standard fusion argument using Lemma~\ref{l:fusion}
shows that in order to prove the theorem, it suffices to prove the
following preservation property:
\begin{condi}\label{c:pres}
Suppose that
\begin{itemize} 
\item  $p \forces{\Poset_\alpha}{{\mathcal U}
\subseteq Z \text{ is  dense open}}$
\item  $\delta $ and $\epsilon$ are real numbers greater than 0
\item  $\Bar{\gamma} : [a,b] \to Z$ is a rectifiable curve
\item  $k\in\omega$
\item  $\Gamma \in [\omega]^{\leq k}$ 
\item  $A$ is a finite subset of $[a,b]$.
\end{itemize}
Then there is some $q \leq_{\Gamma,k}
p$ as well as a rectifiable curve ${\gamma}: [a,b] \to Z$ and a
finite $B\subseteq [a,b]$  such that
\begin{itemize}
\item $q\forces{\Poset_\alpha}{B \cap {\gamma}^{-1}
({\mathcal U}) \cap B \neq\emptyset}$
\item $\lambda_{\mathfrak Z}({\gamma}) < \lambda_{\mathfrak
Z}(\Bar{\gamma}) + \delta$
\item $\norm{\gamma - \Bar{\gamma}} < \epsilon$
\item $\Bar{\gamma}(z) = \gamma(z)$ for each $ z\in A$.
\end{itemize}
\end{condi}
To see that establishing  Condition~\ref{c:pres}  suffices it may,
without loss of generality, be assumed that
$$p \forces{\Poset_\alpha}{{\mathcal U} =
\bigcap_{n\in\omega}{\mathcal U}_n \text{ and each }{\mathcal U}_n
\text{ is dense open}}$$ and that $\{(a_n,b_n,e_n)\}_{n\in\omega}$ is an
enumeration of $\Rationals\times\Rationals\times \omega$.
Suppose that sequences $\{p_n\}_{n\in k}$, $\{\Gamma_n\}_{n\in k}$,
$\{A_n\}_{n\in\omega}$  and
$\{\gamma_n\}_{n\in k}$ have been constructed. 
Let $\Gamma_k\supseteq \Gamma_{k-1}$ be chosen by some scheme that
will ensure that   
$\bigcup_{n\in\omega}\Gamma_n =\bigcup_{n\in\omega}\dom(p_n)$ and
$\card{\Gamma_k} \leq k$ and let $A =
\bigcup_{n\in k}A_k\cap [a_k,b_k]$.
Using Condition~\ref{c:pres} it is possible to construct 
$p_k$, $A_k$ and
$\gamma_k$ such that  
\begin{itemize}
\item if $(a_k,b_k)$ is a nonempty interval in $[a,b]$
then $$p_k\forces{\Poset_\alpha}{{\gamma}_k^{-1}
({\mathcal U}_{e_k}) \cap A_k \neq \emptyset}$$
\item $\lambda_{\mathfrak Z}(\gamma_k) < \lambda_{\mathfrak
Z}({\gamma_{k-1}}) + 2^{-k}$
\item $\norm{\gamma_k - {\gamma}_{k-1}} < \theta_k$
\item $\gamma_k(z) = \gamma_{k-1}(z)$ for each $z\in A$
\end{itemize}
where $\theta_k$ has been chosen so that $\theta_k < 2^{-k}$ and,
using compactness and the fact that $\gamma_{k-1}$ is one-to-one, such
that if $\card{x - y} > 1/k$ and $\{x,y\} 
\subseteq (a_k, b_k)$ then $d(\gamma_{k_1}(x),\gamma_{k_1}(y)) > \theta_k$. 
Hence, $\gamma_\omega = \lim_{n\to\infty}\gamma_n$ is a rectifiable curve and, from
the choice of the $\theta_k$, it follows that $\gamma_\omega$ is also
one-to-one. Moreover, by Lemma~\ref{l:fusion} it follows that there is
some condition $p_\omega$ such that $p_\omega \leq p_n$ for all $n \in
\omega$. Hence $p_\omega\forces{\Poset_\alpha}{{\gamma}_k^{-1} 
({\mathcal U}_{e_k}) \cap A_k \neq \emptyset}$ for each $k \in \omega$.
Since the construction guarantees that $\gamma_\omega\restriction A_k
= \gamma_k\restriction A_k $ it follows that
$p_\omega\forces{\Poset_\alpha}{{\gamma}_\omega^{-1} 
({\mathcal U}_{e_k}) \cap A_k \neq \emptyset}$ for each $k \in \omega$.
Hence 
$p_\omega\forces{\Poset_\alpha}{{\gamma}_\omega^{-1} 
({\mathcal U}_n) \cap [a,b] \text{ is dense in }[a,b]}$ for each $n\in \omega$
and so $p_\omega\forces{\Poset_\alpha({\mathbf X},{\mathfrak
Z})}{{\gamma}_\omega^{-1} 
({\mathcal U}) \cap [a,b] \text{ is dense in }[a,b]}$ .

In order to show that it is possible to satisfy Condition~\ref{c:pres}
it suffices to show that it is possible to satisfy the following
condition: 
\begin{condi}\label{c:pres2}
Suppose that
\begin{itemize} 
\item  $p \forces{\Poset_\alpha}{{\mathcal W}
\subseteq Z \text{ is  dense open}}$
\item  $\delta > 0 $ 
\item  $k\in\omega$
\item  $\Gamma \in [\omega]^{\leq k}$ 
\end{itemize}
Then there is some $q \leq_{\Gamma,k}
p$ as well as a rectifiable curve ${\gamma}: [0,1] \to Z$ and a
finite $B\subseteq [0,1]$ such that
\begin{itemize}
\item $q\forces{\Poset_\alpha}{\Bar{\gamma}^{-1}
({\mathcal W}) \cap B \neq\emptyset}$
\item $\lambda_{\mathfrak Z}({\gamma}) <  \delta$
\end{itemize}
\end{condi}
The reason this suffices it that, given $\Bar{\gamma}:[a,b] \to Z$,
$A$, $k$ and $\epsilon > 0$ as in 
Condition~\ref{c:pres} 
it is possible to choose a point  $w\in [a,b]$ and $\beta > 0$ such that
$\beta < \min(\epsilon, \delta/6)$ and
 such that $B_{{\mathfrak Z}}(\Bar{\gamma}(w),\beta)$ is 
disjoint from $\{\Bar{\gamma}(z) : z\in A\}$. Let $\emptyset \neq
[a',b'] \subseteq 
[a,b]$ be such that $\Bar{\gamma}[a',b'] \subseteq
B_{{\mathfrak Z}}(\Bar{\gamma}(w),\beta)$ and
 define 
 ${\mathfrak Z}' = (Z\cap B_{{\mathfrak
Z}}(\Bar{\gamma}(w),\beta),d)$.
Then apply  
Condition~\ref{c:pres2} to find $\gamma:[0,1] \to B_{{\mathfrak
Z}}(\Bar{\gamma}(w),\beta) $ and $B$ such that 
$\lambda_{{\mathfrak Z}'}(\gamma) < \delta/3$
and note that, since $\gamma$ is a curve and not just a set, it follows that
 $\lambda_{\mathfrak Z}(\gamma) =
\lambda_{{\mathfrak Z}'}(\gamma) < \delta/3$.
It is then an easy matter to
reparametrize $\gamma$ so that it can be assumed to have domain
$[a'',b'']$ where $a''$ and 
$b''$ are chosen so that $a' < a'' < b'' < b'$.
Then use  the hypothesis that $\mathfrak Z$ is reasonably geodesic
 to find curves
$\gamma_{a}:[a',a''] \to B_{{\mathfrak
Z}}(\Bar{\gamma}(w),\beta)$ and $\gamma_{b}[b'',b'] \to B_{{\mathfrak
Z}}(\Bar{\gamma}(w),\beta)$ such that 
\begin{align}
\gamma_{a}(a') = &\Bar{\gamma}(a') &\AND& &
\gamma_{a}(a'')& = {\gamma}(a'')\\
\gamma_{b}(b'') =& {\gamma}(b'') &\AND& &
\gamma_{b}(b') &= \Bar{\gamma}(b')\\
\lambda_{{\mathfrak Z}}(\gamma_{a}) <& \delta/3 &\AND& &
\lambda_{{\mathfrak Z}}(\gamma_{b}) &< \delta/3
\end{align}
Hence $\Bar{\gamma}\cup \gamma\cup \gamma_{a'}\cup \gamma_{b'}$ is the
desired curve and, after reparametrizing, $B$ is the desired finite set.

To verify that Condition~\ref{c:pres2} can be satisfied, let
  $$p \forces{\Poset_\alpha}{{\mathcal W}
\subseteq Z \text{ is  dense open}}$$
  $\delta > 0 $,  
  $\Gamma \in [\omega]^{<\aleph_0}$ 
and $k\in\omega$. Using Corollary~\ref{l:suff} choose $q' \leq_{\Gamma,k}
p$ and a tree of trees $\mathcal S$ which is $k$-sufficient for $q'$
and $\Gamma$. Let $M = \card{\mathcal S}$ and let $\{\sigma_i\}_{i\in M}$
enumerate $\Top({\mathcal S})$.

Using hypothesis~\ref{h:3} choose  $\{x_i\}_{i\in M}$ and  and let
$\theta >0$ be such that 
for any choice of $$(\gamma_1,\gamma_2\dots,\gamma_{M}) \in
\prod_{i\in M}{\mathbf C}(B_{\mathfrak Z}(x_i,\theta))$$ there exists
a single curve 
$\gamma$ such that 
$\lambda_{\mathfrak Z}(\gamma) < \delta$ and $\gamma \supseteq
\bigcup_{i= 1}^M\gamma_i$. 
By shrinking $\theta$ if necessary, it may also be assumed that
$B_{\mathfrak Z}(x_i,\theta)\cap B_{\mathfrak Z}(x_i,\theta)=
\emptyset$.
 Now construct, in $M$ 
steps, conditions $\{q_i\}_{i\in M+1}$  as well as 
finite sets $\{A_i\}_{i\in M+1}$ such that: 
\begin{gather}
\label{ih:1} q_0 = q'\\
\label{ih:9} {\mathcal S}\text{ is good for } q_i \AND \Gamma\\
\label{ih:2} q_{i+1} \sqsubseteq_{{\mathcal S},\sigma_i,M+k} q_i\\
\label{ih:5} q_{i+1}[\Gamma,\sigma_i] \forces{\Poset_\alpha({\mathbf
X},{\mathfrak Z})}{{\mathcal W} \cap A_i \neq 
\emptyset}\\
\label{ih:6} \lambda_{{\mathfrak Z}_i}(A_i) < \theta\\
\label{ih:7} A_i \subseteq { Z}_i
\end{gather}

Before proceeding it is worth pausing to see what would be
accomplished by this construction. It follows from
 Lemma~\ref{l:newkey} that
$q_M\leq_{\Gamma,k} q' $. Using the properties of the family
$\{{\mathcal V}_{i}\}_{i\in M}$ obtained from hypothesis~\ref{h:3} it is
possible to find a curve $\gamma$ such that $\lambda_{\mathfrak Z}(\gamma) < \delta$ and $\cup_{i\in M}A_i \subseteq \gamma$. Let $B =
\gamma^{-1}(\cup_{i\in M}A_i)$ and note that
$q_M\forces{}{B\cap \gamma^{-1}({\mathcal W})) \neq \emptyset}$ because,
by Lemma~\ref{l:maxac}, if $q \leq q_M$ then there is some $i\in M$
such that  $q \leq q_M[\Gamma,\sigma_i] \leq q_{i+1}[\Gamma,\sigma_i]$ and 
$ q_{i+1}[\Gamma,\sigma_i]\forces{}{A_i\cap {\mathcal W} \neq \emptyset}$.

Hence, all that has to be shown is how to construct $q_{i+1}$ given $q_i$.
In order to do this, it suffices to find $\Tilde{q} \leq
q_i[\Gamma,\sigma_i]$ and $A_i\subseteq Z_i$
such that  for each $\alpha \in \Gamma$
\begin{equation}\label{e:c4}
\Tilde{q}\forces{\Poset_\alpha}{(\forall t \in
\Tilde{q}(\alpha)) \text{ either }
\nu^{**}_{\Tilde{q}}(t) \geq M +k\OR
\nu^{**}_{\Tilde{q}}(t) \geq \nu^{**}_{q_i[\Gamma,\sigma_i]}(t) -
1}\end{equation}
and  
$\Tilde{ q} \forces{}{{\mathcal W} \cap A_i \neq
\emptyset}$ and condition~\eqref{ih:6}
is satisfied since, once this has been accomplished, all that remains
to be done is to choose $q_{i+1}$ using
Lemma~\ref{l:nk2} such that that $$q_{i+1}\sqsubseteq_{{\mathcal
S},\sigma_i,M+k} q_i$$ $\mathcal S$ is good for $q_{i+1}$ and
$\Gamma$ and $q_{i+1}[\Gamma,\sigma_i] \leq \Tilde{q}$.

Note that $1\forces{\Poset_{\omega_2}}{\Dot{x} \in V\cap {\mathcal U}
\cap Z_i}$ for 
some name $\Dot{x}$.  Then use Lemma~\ref{l:x} to find
$\Tilde{\Tilde{q}} \leq_{\Gamma,M+k} q_i[\Gamma,\sigma_i]$ and 
a tree of 
trees $\Tilde{\mathcal S}$  which is $(M+k)$-sufficient for both
$\Tilde{\Tilde{q}}$ and $\Gamma$ 
and $q_{i}[\Gamma,\sigma_i]$ and $\Gamma$ such
that $$\Tilde{\Tilde{q}}[\Gamma,\tau] \forces{\Poset_{\omega_2}}{\Dot{x} =
\Check{x}_\tau}$$ for 
each $\tau \in \Top(\Tilde{\mathcal S})$.

Note that 
\begin{equation}\label{e:rr}
\sakne_*(\Tilde{\mathcal S}) \geq \hght_*({\mathcal S})
\geq \Gamma = \hght{\Tilde{\mathcal S}}
\end{equation} and,
 moreover, since $\mathcal S$ is
$k$-sufficient, it follows 
that if $\tau \in \Tilde{\mathcal S}$ and $t
\not\subseteq\sakne(\tsuc_{\Tilde{\mathcal S}}(\tau))$ then
$\nu^{**}_{\Tilde{\mathcal S},\tau}(t) \geq k \geq \Gamma =
\hght{\Tilde{\mathcal S}}$ and hence $\nu_*(\Tilde{S}) >
\hght(\Tilde{S})$. Let the function $\Phi$ be defined on 
$\Top({\Tilde{\mathcal S}})$ by $\Phi(\tau) = x_\tau$.  
It is therefore possible to apply Lemma~\ref{l:m} to $\Tilde{\mathcal
S}$, $\Phi$ and the metric space ${\mathfrak Z}_i = (Z_i,d)$ to find 
 a nonempty tree of trees $\Bar{\mathcal S} \prec_f \Tilde{\mathcal
S}$ such that   
\begin{equation}\label{ih:6f}
\lambda_{{\mathfrak Z}_i}(\{x_\tau\}_{\tau\in\Top({\Bar{\mathcal S}})})) <
\diam({\mathfrak Z}_i) \leq \theta
\end{equation}

Then define $\Tilde{q}$ by
$$\Tilde{q}(\alpha) = \begin{cases}
\Tilde{\Tilde{q}}(\alpha) & \IF \alpha \notin \Gamma\\
\bigcup\{\Tilde{\Tilde{q}}(\alpha)\angbr{s} : s\in \tsuc_{\Bar{\mathcal
S}}(\tau) \AND 
& \IF \alpha \in \Gamma\\
\card{\tau} = \card{(\alpha + 1) \cap \Gamma} \AND 
\Tilde{q}[\Gamma,\tau]\restriction \alpha \in \Dot{G}\cap \Poset_\alpha\} &
		      \end{cases}$$
where $\Dot{G}$ is a name for the generic filter on $\Poset_{\omega_2}$.
It will be shown that $\Tilde{q}$ and $A_i =
\{x_\tau\}_{\tau\in\Top({\Bar{\mathcal S}})}$  are as desired.
Condition~\eqref{ih:6} follows immediately from \eqref{ih:6f}. In
order to check that \eqref{ih:9} holds appeal to Lemmma~\ref{new:tech}
to conclude that $\Bar{\mathcal S} $ is good for $\Tilde{q}$ and $\Gamma$.
Hence $\Tilde{q}\forces{}{{\mathcal W}\cap A_i \neq
\emptyset}$ because, if $ r \leq \Tilde{q}$ then
Lemma~\ref{l:maxac} guarantees that there is
$ \sigma \in \Top(\Bar{\mathcal S})$ such that $r \leq
\Tilde{q}[\Gamma,\sigma] \leq \Tilde{\Tilde{q}}[\Gamma,\sigma]$ and
$\Tilde{\Tilde{q}}[\Gamma,\sigma]\forces{}{\Dot{x} = \Check{x}_\sigma \in
A_i\cap {\mathcal W}}$. 
 Therefore all that has to be checked is that \eqref{e:c4} holds
for all $\alpha \in \Gamma$. 

This follows from the fact that
$\Bar{\mathcal S} \prec_f \Tilde{\mathcal S}$. To
see that this is so notice that, since $\Bar{\mathcal S}$ is good for
$\Tilde{q}$ and $\Gamma$ it suffices to show that for any $\alpha\in
\Gamma$ and any $\tau \in \Bar{\mathcal S}$ such that $\card{\tau} =
\card{\Gamma \cap (\alpha + 1)}$
$$
\Tilde{q}[\Gamma,\tau]\restriction \alpha\forces{\Poset_\alpha}{(\forall t \in
\Tilde{q}(\alpha)) \text{ either }
\nu^{**}_{\Tilde{q}}(t) \geq M +k\OR
\nu^{**}_{\Tilde{q}}(t) \geq \nu^{**}_{q_i[\Gamma,\sigma_i]}(t) -
1}$$ Let such a $\tau$ and $\alpha \in \Gamma$ be fixed. Since
$\Tilde{\mathcal S}$ is $(M+k)$-sufficient $\Tilde{\Tilde{q}}$ and
$\Gamma$ and, furthermore, $\Bar{\mathcal S} \prec
\Tilde{\mathcal S}$ it follows that  
$$\Tilde{q}[\Gamma,\tau]\restriction \alpha\forces{\Poset_\alpha}{(\forall t \in
\Tilde{q} \setminus \tsuc_{\Tilde{\mathcal S}}(\tau) )
\nu_{\Tilde{q}(\alpha)}^{**}(t) =
\nu_{\Tilde{\Tilde{q}}(\alpha)}^{**}(t) \geq M+k}$$ Therefore it
suffices to show that
$$
\Tilde{q}[\Gamma,\tau]\restriction \alpha\forces{\Poset_\alpha}{(\forall t \in
\tsuc_{\Tilde{\mathcal S}}(\tau)) 
\nu^{**}_{\Tilde{q}}(t) \geq \nu^{**}_{q_i[\Gamma,\sigma_i]}(t) -
1}$$ and, since $\Tilde{\mathcal S}$ is good for $q_i[\Gamma,\sigma_i]$ and
$\Gamma$, this is the same as showing that 
$$\nu^{**}_{\Bar{\mathcal S}, \tau}(t) \geq
 \nu^{**}_{\Tilde{\mathcal S}, \tau}(t) -1$$ for each $ t \in
\tsuc_{\Bar{\mathcal S}}(\tau)$. 
For a fixed  $ t \in \tsuc_{\Bar{\mathcal S}}(\tau)$
 the definition of $\Bar{\mathcal S}\prec_f \Tilde{\mathcal
S}$ guarantees that 
\begin{equation}\begin{split}
\nu^{**}_{\Bar{\mathcal S}, \tau}(t) = & \frac{\nu^{*}_{\Bar{\mathcal
S}, \tau}(t)}{(M_{\card{t}}({\mathbf X}))^{{\card{t}}^2}}
\geq 
\frac{\nu^{*}_{\Tilde{\mathcal
S}, \tau}(t) - f(\hght({\Bar{\mathcal S}}),\card{\tau},\card{t})}
{(M_{\card{t}}({\mathbf X}))^{{\card{t}}^2}}\\
= &
\frac{\nu^{*}_{\Tilde{\mathcal
S}, \tau}(t) - M_{\card{t}}({\mathbf X})^{\hght({\Bar{\mathcal S}})
(\hght({\Bar{\mathcal S}}) - \card{\tau})}}
{(M_{\card{t}}({\mathbf X}))^{{\card{t}}^2}}
		\end{split}
\end{equation}
Recalling that $\card{t} \geq  \hght({\mathcal S}) = 
\hght(\Bar{\mathcal S})$ from \eqref{e:rr}
yields that $$\nu^{**}_{\Bar{\mathcal S}, \tau}(t) \geq
\frac{\nu^{*}_{\Tilde{\mathcal
S}, \tau}(t)}
{(M_{\card{t}}({\mathbf X}))^{{\card{t}}^2}} - 1 =
\nu^{**}_{\Tilde{\mathcal S}, \tau}(t) - 1
$$
which establishes \eqref{e:c4}.
\end{proof}

\begin{corol}\label{c:main}
It is consistent with set theory that there are two complete,
separable metric spaces $\mathfrak Y$ and $\mathfrak Z$ such that the
union of
any $\aleph_1$ rectifiable curves in $\mathfrak Y$ is  meagre yet there
are $\aleph_1$ rectifiable curves in $\mathfrak Z$ is whose union is
not  meagre. 
\end{corol}
\begin{proof}
Let ${\mathfrak Y} = \ell^2$ and ${\mathfrak Z} = [0,1]^3$.
From Lemma~\ref{cor:norms} it follows that if $\Poset_{\omega_2}$ is the
countable support iteration of $\Poset({\mathbf X},
{\mathfrak Y}, {\mathfrak Z})$ to $\omega_2$ then, if $G$ is
$\Poset_{\omega_2}$ generic over $V$ then  the
union of
any $\aleph_1$ rectifiable curves in $\mathfrak Y$ is  meagre.
On the other hand, it follows from Theorem~\ref{t:main} that the
union of the ground model rectifiable curves in $\mathfrak Z$ is not meagre.
\end{proof}
\section{Measure and rectifiable curves}
This section contains a result due to Shelah which is somewhat, but
not quite, the measure analogue of Theorem~\ref{t:main}. The measure
analogue would assert the consistency of all sets of reals of size
$\aleph_1$ having measure zero yet there existing a family of
rectifiable curves whose union is not of measure zero.
\begin{theor}[S. Shelah]\label{t:shel}
It is consistent with set theory that $\non(Null) > \aleph_1$ yet
$\non({\mathcal R}(\Reals^2)) = \aleph_1$.
\end{theor}
\begin{proof}
For certain pairs of functions $(g,g^*) \in \fomom\times\fomom$ a
partial order $\Sacks_{g,g^*} $ is defined on pages 348-349 of
\cite{ba.ju.book}. It is shown that $\Sacks_{g,g^*} $ is proper and
$(g^*,g)$-bounding\footnote{The careful reader will notice that, in
fact, Lemma~7.3.23 claims that  $\Sacks_{g,g^*} $ is
$(g,g^*)$-bounding but this clearly a misprint.} and hence, from
Lemma~7.2.19, which establishes that the countable support iteration
of proper $(f,h)$-bounding partial orders is $(f,h)$-bounding, it
follows that the countable support iteration of length $\omega_2$ of
$\Sacks_{g,g^*} $ is also $(g^*,g) $-bounding. Furthermore, it follows
directly from Theorem~7.3.21 that if $G$ is generic for this
iteration over a model $V$ then $V[G]$ is a model of $\non(null) >
\aleph_1$.

Hence it suffices to show that if $g$ and $g^*$ are chosen
appropriately then $V[G]$ will also be a model of $\non({\mathcal
R}(\Reals^2)) = \aleph_1$. In particular, it will be shown that if
$\{\gamma_i\}_{i\in\omega}$ is any family of rectifiable curves in
$\Reals^2$ in the 
generic extension then there is a point $x\in V\cap \Reals^2$ such
that $x \notin \bigcup_{i\in\omega}\gamma_i$. Let $f$, $g$ and $g^*$ be
chosen so that $g$ and $g^*$ satisfy (1) and (2) on page 348 of
\cite{ba.ju.book} and, in addition, so that $$g^*(n) =
\card{[f(n)^2]^{<5f(n)}}$$
$f(0) = 5$, $g(0) = 1$
and $$\frac{f(n+1)}{f(n)^2} > 5g(n+1)$$
as well. To see that this works,  let $\{\gamma_i\}_{i\in\omega}$ be
a family of rectifiable curves in the generic extension and, without
loss of generality, suppose that none of the curves has length greater
than 1. Attention will be focussed on the unit square $[0,1]^2$.

Before continuing, observe that if the unit square  is
covered by $m^2$ congruent closed squares with sides of length $1/m$
and $\gamma$ is 
any rectifiable curve of length  less than $1$ then there are at most
${5m}$ squares which the curve intersects. The way to see this is to choose
a point from the curve from each square which the curve intersects.
Let these points be $\{\gamma(x_i)\}_{i=0}^k$ where $x_i \leq
x_{i+1}$. Since no 5 of the squares can have a common point of
intersection, it follows that for any $i\in k-5$ there must be some
$\{i,j\} \in [5]^2$ such that the  ${\lVert {\gamma(x_i) -
\gamma(x_j)} \rVert} \geq 1/m$. Since $\sum_{i\in k}{\lVert {\gamma(x_i) -
\gamma(x_{i+1})} \rVert} < 1$ it follows that $k < 5m$.

Now let ${\mathcal L}(n)$ be the natural  family of size $f(n)^2$
consisting of congruent closed squares of length $1/f(n)$  whose union
covers the unit square. For each $n \in \omega$ let $S_n$ be the set
of squares in ${\mathcal L}(n)$ which the curve $\gamma_n$ intersects.
Since $S_n \in [{\mathcal L}(n)]^{<5f(n)}$ it follows that there is
$H: \omega \to [[{\mathcal L}(n)]^{<5f(n)}]^{g(n)}$ in $V$ such that $S_n \in 
H(n)$ for all $n\in\omega$. Let $S^*_n = \cup H(n)$. It follows that
$\card{S_0^*} < 5\cdot 5$ and so there is some square $Q_0 \in {\mathcal
L}(n)$ such that $Q_0\notin S^*_0$. Suppose that
inductively square $Q_0 \supseteq Q_1 \supseteq \ldots \supseteq Q_n$
have been chosen so that $Q_i\notin S^*_i$. The number of squares in
${\mathcal L}(n+1)$ which are subsets of $Q_n$ is
$\frac{f(n+1)^2}{f(n)^2}$ and $\card{S^*_{n+1}} \leq 5f(n+1)g(n+1)$.
Since $$\frac{f(n+1)}{f(n)^2} > 5g(n+1)$$ it follows that it is
possible to find $Q_{n+1}\subseteq Q_n$ which does not belong to
$S^*_{n+1}$. Letting $x = \bigcap_{n\in\omega}Q_n$ yields that
$x\notin \bigcup_{n\in\omega}S^*_n$ and hence $x\notin
\bigcup_{n\in\omega}\gamma_n$ 
\end{proof}
\section{Open questions}
The proof of Theorem~\ref{t:main} relied heavily on the hypothesis
that $\mathfrak Z$ is reasonably geodesic and this, in turn, implies
some form of dimensionality  greater than 3. This raises the obvious
question of whether this is only a technical problem or if there is
some intrinsic geometric reason for the hypothesis. Answering the
following questions may shed some light on this.
\begin{quest}
Is it consistent that $\non(Meagre) > \aleph_1$ yet there is a family
of $\aleph_1$  rectifiable curves in $\Reals^2$ whose union is not
meagre?
\end{quest}
\begin{quest}
Is it consistent that that the union of any $\aleph_1$ rectifiable
curves in $\ell^2$ is meagre  yet there is a family
of $\aleph_1$  rectifiable curves in $\Reals^2$ whose union is not
meagre?
\end{quest}
\begin{quest}
Are there integers $n$ and $k$ such that it consistent that that the
union of any $\aleph_1$ rectifiable 
curves in $\Reals^n$ is meagre  yet there is a family
of $\aleph_1$  rectifiable curves in $\Reals^k$ whose union is not
meagre?
\end{quest}
At this point it should be pointed out that there is an obvious
strategy for proving some monotonicty results by taking a family of
rectifiable curves in $\Reals^{n+1}$ whose unions is not meagre and
projecting them down to $\Reals^n$ the curves will still be
rectifiable and their union will not be meagre. However the projected
curves may not be one-to-one.  It would be interesting to see if
theorem~\ref{t:main} can be improved so that, in the case ${\mathfrak
Z} = \Reals^n$, the rectifiable curve
mentioned in the conclusion is nicely embedded. There are various
interpretations what nicely embedded might mean here. The best would
be a linear embedding of the unit interval since this would solve
Komjath's question mentioned in the introduction. A less ambitious
alternative is the following.
\begin{quest}
Is it consistent that $\non(Meagre) > \aleph_1$ yet there is a family
of $\aleph_1$  rectifiable curves in $\Reals^n$ whose union is not
meagre and all of which are homotopic to the unit interval?
\end{quest}
In this same spirit, one can ask for smoothness properties.
\begin{quest}
Is it consistent that $\non(Meagre) > \aleph_1$ yet there is a family
of $\aleph_1$  differentiable curves in $\Reals^n$ whose union is not
meagre.?
\end{quest}
Of course, one can ask for higher orders of smoothness as well. Nothing 
is known in any of these cases either.

Finally, it should be mentioned that al of the questions considered
here have their measure analogues. Indeed, rectifiable curves are
special examples of sets of Hausdorff dimension 1 in $\Reals^2$. One
can oslo consider the $\sigma$-ideal generated by Borels sets of
Hausdorff dimension $p$ in $\Reals^n$ and ask similar questions in
this a context. For some progress in this direction as well as further
refernces see \cite{step.37}.

\end{document}